%% file: main.tex
\expandafter\def\csname ver@fixltx2e.sty\endcsname{}
\documentclass[a4paper,oneside,preprint,12p]{elsarticle}

\usepackage{url}
\usepackage{microtype}

\usepackage{dblfloatfix}
\usepackage[utf8]{inputenc}
\usepackage{natbib}
\usepackage{amsmath}
\usepackage{amsthm}
\usepackage{amssymb}
\usepackage{amsfonts}
\usepackage{bm}
\usepackage{siunitx}
\usepackage{epstopdf}          
\usepackage{flushend}
\usepackage{parskip}
\usepackage{hyperref}
\usepackage[hyperref,dvipsnames]{xcolor}

\definecolor{bluePolimi}{RGB}{22, 44, 80}
\definecolor{lightBluePolimi}{RGB}{91, 122, 172}
\definecolor{redPolimi}{RGB}{180, 0, 0}
\definecolor{greenPolimi}{RGB}{78, 172, 91}
\definecolor{green2}{RGB}{0, 110, 0}

\makeatletter
\hypersetup{colorlinks=true}
\AtBeginDocument{\@ifpackageloaded{hyperref}
  {\def\@linkcolor{blue}
   \def\@anchorcolor{red}
   \def\@citecolor{red}
   \def\@filecolor{red}
   \def\@urlcolor{redPolimi}
   \def\@menucolor{red}
   \def\@pagecolor{cyan}
\begingroup
  \@makeother\`%
  \@makeother\=%
  \edef\x{%
    \edef\noexpand\x{%
      \endgroup
      \noexpand\toks@{%
        \catcode 96=\noexpand\the\catcode`\noexpand\`\relax
        \catcode 61=\noexpand\the\catcode`\noexpand\=\relax
      }%
    }%
    \noexpand\x
  }%
\x
\@makeother\`
\@makeother\=
}{}}
\makeatother

\usepackage[pagewise]{lineno}
\usepackage{framed} 
\usepackage{multicol} 
\usepackage{tabularx} 
\usepackage{longtable} 
\usepackage[ruled,vlined]{algorithm2e}

\journal{Applied Mathematical Modelling}
\usepackage{amssymb}

\usepackage[flushright]{rotating}

\usepackage{subcaption} 
\usepackage[printonlyused]{acronym} 
\usepackage{mathtools} 

\usepackage{tikz}
\usetikzlibrary{fit}
\usetikzlibrary{positioning}
\usetikzlibrary{shapes,arrows, arrows.meta}

\DeclareMathOperator*{\argmax}{arg\,max}
\DeclareMathOperator*{\argmin}{arg\,min}

\DeclareMathOperator*{\arginf}{arg\,inf}
 \newcolumntype{M}[1]{>{\centering\arraybackslash}m{#1}}
 \usepackage{booktabs} 
 \usepackage{comment}

\renewcommand{\vec}[1]{\mathbf{#1}}

\newcommand{\module}[1]{\left|{#1}\right|}
\newcommand{\norma}[1]{\left|\left|{#1}\right|\right|}

\newcommand{\scalarprod}[2]{\left({#1},\,{#2}\right)}

\newcommand{\dpart}[2]{\frac{\partial {#1}}{\partial {#2}}}

\begin{document}

\newtheorem{theo}{Theorem}
\theoremstyle{definition}
\newtheorem{obs}{Remark}
\newtheorem{Def}{Definition} 

\begin{frontmatter}

\title{Multi-Physics Model Bias Correction with Data-Driven Reduced Order Modelling Techniques: \\Application to Nuclear Case Studies}

\author[First]{Stefano Riva}
\author[First]{Carolina Introini}
\author[First, cor1]{Antonio Cammi}

\cortext[cor1]{Corresponding author. Email address: antonio.cammi@polimi.it}

\address[First]{Politecnico di Milano, Dept. of Energy, CeSNEF-Nuclear Engineering Division, Nuclear Reactors Group - via La Masa, 34 20156 Milano, Italy}

\begin{abstract}
Nowadays, interest in combining mathematical knowledge about phenomena and data from the physical system is growing. Past research was devoted to developing so-called high-fidelity models, intending to make them able to catch most of the physical phenomena occurring in the system. Nevertheless, models will always be affected by uncertainties related, for example, to the parameters and inevitably limited by the underlying simplifying hypotheses on, for example, geometry and mathematical equations; thus, in a way, there exists an upper threshold of model performance. Now, research in many engineering sectors also focuses on the so-called data-driven modelling, which aims at extracting information from available data to combine it with the mathematical model. Focusing on the nuclear field, interest in this approach is also related to the Multi-Physics modelling of nuclear reactors. Due to the multiple physics involved and their mutual and complex interactions, developing accurate and stable models both from the physical and numerical point of view remains a challenging task despite the advancements in computational hardware and software, and combining the available mathematical model with data can further improve the performance and the accuracy of the former.

This work investigates this aspect by applying two Data-Driven Reduced Order Modelling (DDROM) techniques, the Generalised Empirical Interpolation Method and the Parametrised-Background Data-Weak formulation, to literature benchmark nuclear case studies. The main goal of this work is to assess the possibility of using data to perform model bias correction, that is, verifying the reliability of DDROM approaches in improving the model performance and accuracy through the information provided by the data. The obtained numerical results are promising, foreseeing further investigation of the DDROM approach to nuclear industrial cases.
\end{abstract}


\end{frontmatter}

\section{Introduction}

The issue of state estimation is a crucial point in the safety, control and monitoring of engineering systems. Mathematical models can retrieve the complete spatial dependence of the variable of interest; however, they are affected by the modelling assumptions related to geometry and mathematical equations and the uncertainty on the values of the parameters characterising the system. The typical representation of these models is through a system of Partial Differential Equations (PDEs) \cite{quarteroni_numerical_2016}, which requires powerful computational hardware to provide an accurate numerical solution: this model is typically called \textit{high-fidelity} of Full Order Model (FOM), and the required computational resources and the time needed to retrieve a solution represent a prominent limitation to their application in multi-query and real-time scenarios \cite{rozza_model_2020}. Another possibility for monitoring engineering systems is through data, i.e., the direct evaluation of the variables of interest; however, data are typically affected by random or systematic noise, and they cannot give complete coverage of the spatial domain \cite{argaud_sensor_2018}. Thus, in the last few years, a hybrid approach based on Data Assimilation (DA), the \textit{modelling from data}, has reached widespread interest \cite{brunton_data-driven_2022}.

Data Assimilation \cite{carrassi_data_2018} methods allow for integrating these two sources of information, the mathematical model and the data, exploiting the advantages of both approaches whilst minimising their drawbacks. Traditional DA techniques build optimisation problems and try to minimise the distance between the mathematical \textit{background} model and the experimental data; this approach requires several evaluations of the former, which means solving the system of PDEs for the numerical optimisation several times, thus making this framework unfeasible, for example, for multi-query scenarios and optimisation applications. Reduced Order Modelling (ROM) \cite{rozza_model_2020, lassila_model_2014, quarteroni2015reduced} represents an attractive solution to reduce the computational complexity of such problems whilst keeping the accuracy to the desired level by substituting the ROM with a surrogate low-dimension model with much lower complexity. Two major categories of ROM methods exist: intrusive algorithms use a Galerkin projection and require knowledge of the governing equations \cite{lorenzi_pod-galerkin_2016, lorenzi_reduced_2017, stabile_pod-galerkin_2017}; non-intrusive methods retrieve the state estimation starting from available data and by solving either a linear system \cite{gong_empirical_2016, maday_general_2008, maday_generalized_2013} or by using Machine Learning \cite{ortali_gaussian_2022, gong_data-enabled_2022} techniques or Neural Networks \cite{hesthaven_non-intrusive_2018, wang_non-intrusive_2019, chen_physics-informed_2021}. This latter approach is more suited for integrating measurements to a reduced order model in the so-called Data-Driven Reduced Order Modelling framework (DDROM), even though examples of hybrid intrusive approaches exist \cite{introini_KalmanFilter, Introini_PhDTesi, He_KalmanROM}. 

\begin{figure}[htbp]
    \centering
    \includegraphics[width=1.\linewidth]{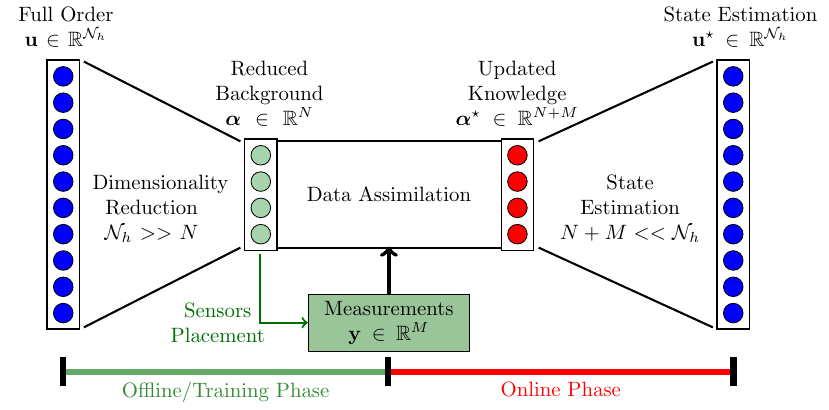}
    \caption{Conceptual scheme of Data-Driven Reduced Order Modelling.}
    \label{fig: DA-ROM}
\end{figure}

Moreover, non-intrusive strategies allow the inclusion of algorithms for sensor placement \cite{argaud_sensor_2018, BinevMula2018_sensors, binev_data_2017, aretz_greedy_2024}. Figure \ref{fig: DA-ROM} reports the idea of DDROM with sensor placement: in the offline (training) phase, a dimensionality reduction process retrieves a reduced coordinate system onto which encodes the information of the mathematical model; the sensor positioning algorithm then uses this set to select the optimal location of sensors according to some optimality criterion, which depends on the adopted algorithm. In the online phase, the DA process begins, retrieving a novel set of reduced variables and then computing the reconstructed state through a decoding step.

Within these hybrid DDROM approaches, two important examples are the Generalised Empirical Interpolation Method (GEIM) \cite{maday_generalized_2013, maday_generalized_2015, maday_convergence_2016, gong_generalized_2022} and the Parameterised-Background Data-Weak (PBDW) \cite{maday_parameterized-background_2014, maday_yvon_pbdw_2015, maday_adaptive_2019, taddei_model_2016, maday_pbdw_2015, haik_real-time_2023}; in particular, the latter comes with a complete theoretical description of the logic behind DDROM for integration of un-modelled physics using measurements and model bias correction. Literature examples of GEIM and PBDW applied to complex engineering systems such as nuclear reactors \cite{gong_generalized_2022, Introini_PhDTesi, cammi_indirect_2023} and industrial experimental facilities \cite{taddei_model_2016, haik_real-time_2023, riva_hybrid_2023-1} exist, proving their efficiency and their reliability. 

This work focuses on applying these techniques in the context of Multi-Physics (MP) modelling of nuclear reactors \cite{demaziere_6_2020}, which presents several coupling effects between different physical phenomena which play a central role in the design and control of the system itself: the most significant example of such coupling is the one between neutron physics and Thermal-Hydraulics (TH) \cite{DuderstadtHamilton}. Indeed, the power generated by fission events induces a temperature evolution which affects the fluid dynamics, especially in Circulating Fuel Reactors \cite{Aufiero2014}, through buoyancy effects and changes in the material properties for the interaction with neutrons \cite{hebert_applied_2009}. Until a few years ago, mathematical modelling of nuclear reactors focused on accurate numerical codes describing a single physics (SP), such as OpenFOAM \cite{OPENFOAM} for TH and Serpent \cite{Serpent} for neutronics, with the coupling information given as boundary condition. With the advancements in computational hardware, the MP approach is overcoming the SP modelling strategy. MP models can achieve coupling between physics in two ways: either through the development of interfaces between single-physics codes \cite{SerpentOpenFOAM_Godiva, Castagna_Pin_SerpentOpenFOAM, SerpentOpenFOAM_PB} or by gathering all the information inside a single new environment \cite{fiorina_gen-foam_2015}. The latter path is to become state-of-the-art in the nuclear field, but using this approach means losing most knowledge and validation studies performed on SP codes; in addition, the computational cost for such complex models remains high. In this context, DDROM methods can introduce un-modelled physics and unanticipated uncertainty in the existing and already-validated SP codes through data \cite{maday_parameterized-background_2014} through the data: in particular, validated SP codes generate the training snapshots to obtain an SP representation in a reduced coordinate system \cite{quarteroni2015reduced, brunton_data-driven_2022}, introducing the coupling information through 'boundary conditions'; then, some MP data (which can come either from experiments or some simpler MP model) update and correct the model bias. In the literature, other approaches can be adopted to combine different levels of fidelity (e.g, the SP and the MP) based on Gaussian Process Regression \cite{kast_non-intrusive_2020} and Long short-term memory networks \cite{conti_multi-fidelity_2023}, however in this paper it has been assumed that the knowledge of the true MP model is accessible only through data in the form of measurements, making GEIM and PBDW more suited. Therfore, this work follows this concept and aims at assessing their capabilities for online model bias correction for some benchmark MP cases, such as the IAEA 2D PWR \cite{argonne_book, Milonga_code} and TWIGL2D \cite{twigl-Hageman, twigl-Yasinsky}.

The paper is organised as follows: Section \ref{sec-ROM} describes the data-driven ROM methods adopted in this work, with a brief theoretical analysis of the proposed approach for MP modelling of nuclear reactors; Section \ref{sec: MP-case-study} discusses the governing equations and the case studies, whereas Section \ref{sec: num-res} is devoted to the numerical results, focusing on the model bias correction played by the data-driven ROM methods. In the end, in Section \ref{sec: conclusions}, the main conclusions for the paper are drawn.

\section{Reduced Order Modelling Techniques for Data Assimilation}\label{sec-ROM}

This section presents the methods employed in this work, i.e., the Generalised Empirical Interpolation Method (GEIM) and the Parameterised-Background Data-Weak (PBDW) formulation, focusing on their capabilities of correcting model bias. As many ROM techniques, both require a collection of snapshots $\{u(\vec{x};\,\boldsymbol{\mu})\}_{\boldsymbol{\mu}\in \mathcal{D}}$, typically solution of a parameterised PDE, dependent on $\boldsymbol{\mu}\in \mathcal{D}\subset\mathbb{R}^p$, $p\geq 1$, belonging to a suitable functional space $\mathcal{U}\subset L^2(\Omega)$ over a spatial domain $\Omega$ endowed with a scalar product, to generate a set of basis functions and basis sensors: the former are necessary to get a reduced representation of the solution space, whereas the latter allows integrating data/measures for the model bias correction. This selection occurs in the offline phase, which is computationally expensive but performed only once. During the online phase, instead, the data $\vec{y}\in\mathbb{R}^M$ are collected from a physical system to retrieve a full-state estimation by solving a linear system of reduced dimensions (compared to the FOM). Thus, the online phase ensures fast yet accurate enough simulations suitable for real-time and multi-query scenarios. The separation between a computationally expensive offline (training) step, performed only once, and a computationally cheap online one is the typical framework of Data-Driven Reduced Order Modelling techniques, also called Hybrid Data Assimilation (HDA) \cite{riva_hybrid_2023}, which allows combining the mathematical (\textit{background}) models and available data.

\subsection{Generalised Empirical Interpolation Method}\label{sec-geim}

The first work discussing the Generalised Empirical Interpolation Method (GEIM) was \cite{maday_generalized_2013}, followed by \cite{maday_generalized_2015, maday_convergence_2016, argaud_sensor_2018} and an extension to noisy data in \cite{argaud_stabilization_2017, gong_generalized_2022, introini_stabilization_2023}. Briefly, GEIM works with a greedy procedure that generates the reduced spaces $X_M$ of finite dimension $M$ in a hierarchical\footnote{The following property holds
\begin{equation*}
    X_1\subset X_2\subset \dots \subset X_{M_{\text{max}}}\subset \mathcal{U}
\end{equation*}
which means that at each step, the space is enriched with new information.}  manner \cite{maday_generalized_2015}:the basis functions $\{q_m(\vec{x})\}_{m=1}^M$, spanning this space $X_M$ are known in literature as \textbf{magic functions} and they contain the spatial behaviour of the snapshots. Hence, any function $u(\vec{x};\,\boldsymbol{\mu})$ can be approximated by an interpolant $\mathcal{I}_M$ 
\begin{equation}
    u(\vec{x};\,\boldsymbol{\mu}) \simeq \mathcal{I}_M[u](\vec{x};\,\boldsymbol{\mu}) = \sum_{m=1}^M \beta_m(\boldsymbol{\mu})\cdot q_m(\vec{x}),
    \label{eqn: sec2-geim-interpolant-def}
\end{equation}
in which the the coefficients $\{\beta_m(\vec{x})\}_{m=1}^M$,  obtained through the solution of an interpolation problem, gives the parametric dependence. Each magic function is associated with a \textbf{magic sensor}, a linear functional $v_m(\cdot;\,\vec{x}_m,\,s):\mathcal{U}\rightarrow \mathbb{R}$ which mathematically represents the data acquisition through sensors. These functionals are generally defined as follows \cite{BinevMula2018_sensors, haik_real-time_2023}:
\begin{equation}
    v_m(u(\vec{x});\,\vec{x}_m,\,s)=\int_\Omega u(\vec{x})\cdot \mathcal{K}\left(\|\vec{x}-\vec{x}_m\|_2, s\right)\,d\Omega,
    \label{eqn: sens-def}
\end{equation}
given $\mathcal{K}\left(\|\vec{x}-\vec{x}_m\|_2, s\right)$ a kernel function which depends on the distance from the centre $\vec{x}_m$ (position of the physical sensor) and the point spread $s$. One of the most common choices for the measurements of scalar fields is the Gaussian kernel normalised with respect to the $L^1$-norm\footnote{This norm is defined for $u\in L^1(\Omega)$ as
\begin{equation*}
    \norma{u}_{L^1(\Omega)} = \int_\Omega |u|\,d\Omega  
\end{equation*}}
\cite{maday_generalized_2015}. These functionals belong to a library $\Upsilon\subset\mathcal{U}'$, given $\mathcal{U}'$ the dual space of $\mathcal{U}$. The greedy offline procedure selects both the magic functions and the magic sensors. 

\subsubsection{Offline Greedy Algorithm}

The greedy algorithm requires as input a set of training snapshots, coming from the \textit{high-fidelity} solution of the parameterised PDEs, $\{u^{(i)} = u(\vec{x}, \boldsymbol{\mu}_i)\}_{i=1}^{N_s}$, with $\boldsymbol{\mu}_i\in\Xi_{\text{train}}$ and $\text{dim}\left(\Xi_{\text{train}}\right)=N_s$: at each step $m$, the procedure selects the magic function and associated sensor to minimise the interpolation error made by $\mathcal{I}_m$. Algorithm \ref{Algo: GEIM-Offline} summarises the greedy procedure.

\begin{algorithm}[htbp] 
	\textbf{Input}\\
	$\quad$ Maximum number of iterations $M_{max}$\;
	$\quad$ Tolerance $\delta$\;
	$\quad$ Training Snapshots $\{u^{(i)} = u(\vec{x}, \boldsymbol{\mu}_i)\}_{i=1}^{N_s}$ with $\boldsymbol{\mu}_i\in\Xi_{\text{train}}$\;
	\textbf{Output}\\
	$\quad$ Magic functions $\{q_{1}(\mathbf{x}),...,q_{M}(\mathbf{x})\}$\;
	$\quad$ Magic sensors $\{v_{1},...,v_{M}\}$\; $ $\\
	\textbf{Initialisation}\\
	$\quad M = 1$, $ E_1 = \delta+1$\;
	\textbf{First iteration}\\
	$\quad u_{1} = \argmax\limits_{i = 1, \dots{N_s}}||u^{(i)}||_{L^2(\Omega)}$\;
	$\quad  v_1 = \argmax\limits_{v^{(k)}\in\Upsilon}\left|v^{(k)}(u_1(\vec{x});\vec{x}_k,s)\right|$\;
	$\quad q_1(\vec{x}) = \frac{u_1(\vec{x})}{v_1(u_1(\vec{x}))}$\;
	\While{($M < M_{max} \quad \& \quad \delta \geq E_{M-1}$)}{
		$\quad M=M+1$\;
		$\quad u_M(\vec{x}) = \argmax\limits_{i = 1, \dots{N_s}}\norma{u^{(i)}-\mathcal{I}_{M-1}[u^{(i)}]}_{L^2(\Omega)}$\;
		$\quad$given $I_{M-1}[u^{(i)}] = \sum_{m=1}^{M-1}\beta_m(\boldsymbol{\mu})\cdot q_m(\vec{x})$\;
		$\quad$and $\sum_{m=1}^{M-1} \beta_m\cdot v_{m'}(q_m) = v_{m'}(u^{(i)})$ with $m'=1, \dots, M-1$\;
		$\quad v_M = \argmax\limits_{v^{(k)}\in\Upsilon}\left|v^{(k)}(u_M-\mathcal{I}_{M-1}[u_M];\vec{x}_k,s)\right|$\;
		$\quad q_M(\vec{x}) = \displaystyle \frac{u_M(\vec{x})-\mathcal{I}_{M-1}[u_M](\vec{x})}{v_M(u_M-\mathcal{I}_{M-1}[u_M])}$\;
		$\quad E_{M} =\argmax\limits_{i = 1, \dots{N_s}}\norma{u^{(i)}-\mathcal{I}_{M}[u^{(i)}]}_{L^2(\Omega)}$\;
	}
	\caption{GEIM greedy procedure - Offline Stage}
	\label{Algo: GEIM-Offline}
\end{algorithm} 
At each time step $m$, the procedure enriches the reduced space by adding the worst reconstructed snapshot, labelled by a 'magic' parameter $\boldsymbol{\mu}_m$; then, retrieving the magic functions forming this space. Moreover, the interpolant $\mathcal{I}_M$ must satisfy the interpolation condition with the selected snapshots $\{u(\vec{x}, \boldsymbol{\mu}_k)\in\mathcal{U}\}_{k=1}^{M}$:
\begin{equation}
    \sum_{m=1}^{M}v_{m'}(q_m)\cdot \beta_m=
    \sum_{m=1}^{M}\mathbb{B}_{m'm}\cdot \beta_m = v_{m'}(u(\vec{x}, \boldsymbol{\mu}_k))\qquad \forall m'=1, \dots, M,
\end{equation}
in which the matrix $\mathbb{B}\in\mathbb{R}^{M\times M}$ is lower triangular, ensuring the well-posedness of the linear system \cite{maday_generalized_2015, maday_convergence_2016}.

\subsubsection{Online Estimation}
In the ROM framework, the online phase in the ROM framework needs to return quickly and efficiently the full-order state estimation of the quantity of interest \cite{quarteroni2015reduced}; thus, it must have a small computational cost. For the GEIM, the full-order state estimation is retrieved given some evaluations of the physical system:
\begin{equation}
    y_m = v_m (u^{\text{true}})+\epsilon_m \qquad m = 1, \dots, M,
    \label{eqn: sec2-online-measure}
\end{equation}
with $\epsilon_m$ representing external disturbances (such as noise). Then, the coefficients $\boldsymbol{\beta}\in\mathbb{R}^M$ are obtained as the solution of the interpolation condition:
\begin{equation}
    \sum_{m'=1}^{M}\mathbb{B}_{mm'}\cdot \beta_{m'} = y_m\qquad \forall m=1, \dots, M.
    \label{eqn: sec2-geim-online}
\end{equation}
Assuming clean data (i.e., $\epsilon_m = 0$), the matrix $\mathbb{B}$ is lower-triangular, thus ensuring the well-posedness of the GEIM procedure  \cite{maday_convergence_2016}, meaning that the interpolant exists and is unique. By removing the assumption of clean data, the error bound increases as more sensors are added\footnote{This trend has also been proved numerically in \cite{argaud_stabilization_2017, introini_stabilization_2023}.}, because in its standard formulation, the GEIM is not robust against noise \cite{argaud_stabilization_2017} and thus it requires stabilisation techniques. This work uses the stabilisation proposed by the authors in \cite{introini_stabilization_2023}, based on the Tikhonov regularisation \cite{Tikhonov_Reg}, which works by weakening the interpolation condition \eqref{eqn: sec2-geim-online} into a least-squares formulation with a penalisation term. This new problem is equivalent to the following:
\begin{equation}
    \left(\mathbb{B}^T\mathbb{B}+\lambda \mathbb{T}^T\mathbb{T}\right)\boldsymbol{\beta} = \mathbb{B}^T\vec{y}+\lambda \mathbb{T}^T\mathbb{T} \overline{\boldsymbol{\beta}},
    \label{eqn: sec2-trgeim-online}
\end{equation}
where $\lambda\in\mathbb{R}^+$ is a regularisation parameter to be suitably calibrated\footnote{If $\epsilon_m$ are uncorrelated Gaussian random variables with variance $\sigma^2$, the optimal value of $\lambda$ is $\sigma$ itself when there are no constraints on the sensors position \cite{introini_stabilization_2023, cammi_icapp2023_esteso}. Otherwise, hyperparameter tuning is required.}, $\overline{\boldsymbol{\beta}}$ is the sample mean of the coefficients of the train set ($\Xi_{\text{train}}$), $\mathbb{T}\in\mathbb{R}^{M\times M}$ is the regularisation matrix, defined as in \cite{Tikhonov_Reg, riva_hybrid_2023}:
\begin{equation}
\begin{split}
    \mathbb{T}_{ij} =\left\{
	\begin{array}{cc}
        |\sigma_{\beta_i}|^{-1} & \mbox{ if } i=j\\ 
		0 & \mbox{ if } i\neq j\\
	\end{array}
	\right.\qquad i,j = 1, \dots, M.
	  \end{split}
\end{equation}
with $\sigma_{\beta_i}$ as the sample standard-deviation of the training coefficients. This version of the method is hereafter called TR-GEIM. Algorithm \ref{Algo: GEIM-Online} summarises the overall online procedure.

\begin{algorithm}[htbp] 
	\textbf{Input}\\
	$\quad$ Magic functions $\{q_{1}(\mathbf{x}),...,q_{M}(\mathbf{x})\}$\;
	$\quad$ Magic sensors $\{v_{1},...,v_{M}\}$\;
	\textbf{Output}\\
	$\quad$ State Estimation $\mathcal{I}_M[u^{true}](\vec{x})$\; $ $\\
	\textbf{Online Estimation}\\
	$\quad$ Acquisition of Experimental Data\\
	$\qquad \left\{y_m = v_m(u^{\text{true}})+\epsilon_m\right\}_{m=1}^M$\;
     \eIf{Noise is False}
    {
        Solve for $\boldsymbol{\beta}$ from the interpolation problem \eqref{eqn: sec2-geim-online}
    }{
        Solve for $\boldsymbol{\beta}$ from the regularised problem \eqref{eqn: sec2-trgeim-online}
    }
	$\quad$Estimate the state of the system\\
	$\qquad \displaystyle u^{\text{true}} \simeq \mathcal{I}_M[u^{\text{true}}](\vec{x}) = \sum_{m=1}^M \beta_m\cdot q_m(\vec{x})$
	\caption{(TR-)GEIM - Online Stage}
	\label{Algo: GEIM-Online}
\end{algorithm} 

Both the offline and online phases are implemented in the Python package \href{https://github.com/ERMETE-Lab/ROSE-pyforce}{pyforce}, built upon the \href{https://fenicsproject.org}{FEniCSx} library \cite{BarattaEtal2023, ScroggsEtal2022, BasixJoss, AlnaesEtal2014}.

\subsection{Parameterised-Background Data-Weak formulation}

The state estimation resulting from the combination between the knowledge of the \textbf{background} mathematical model and experimental data from the system, can be written in a variational data assimilation framework \cite{lorenc_analysis_1986, rabier_variational_2003, carrassi_data_2018}:
\begin{equation}
	u^\star = \argmin\limits_{u\in\mathcal{U}}\xi \norma{u-u^{bk}}_{L^2}^2+\frac{1}{M}\sum_{m=1}^M\left(v_m(u)-y_m^{obs}\right)^2,
	\label{eqn: DA-variational}
\end{equation}
where $u^{bk}$ is the solution of the \textit{background} model and $\xi$ is called the regularising parameter, which weights the relative importance of the background (mathematical model) compared to the experimental data: in particular, its optimal value depends on the model error and the noise level. For a perfect model, high values of $\xi$ are suggested; for perfect measurements, $\xi\rightarrow 0^+$ \cite{maday_adaptive_2019}. Moreover, this optimal value increases with noise and decreases as the best-fit error increases \cite{taddei_model_2016}. The optimisation problem \eqref{eqn: DA-variational} looks for an estimation close enough to the mathematical model and enriched with the novel information of the data.

Since solving the optimisation problem directly is not feasible, ROM techniques are crucial to evaluate quickly the background model \cite{aretz-nellesen_3d-var_2019}. The Parameterised-Background Data-Weak formulation, proposed by \cite{maday_parameterized-background_2014} and then extended in \cite{maday_pbdw_2015, taddei_model_2016, maday_adaptive_2019, haik_real-time_2023}, is a general framework for Hybrid Data Assimilation methods, combining ROM and real data \cite{riva_hybrid_2023}. The solution space of the mathematical model is replaced by a reduced one $Z_N$ spanned by some basis functions $\{\zeta_n\}_{n=1}^N$ so that the state estimation can be decomposed into a background term $z_N$ and an update term $\eta_M$, i.e.
\begin{equation}
\begin{split}
    u^\star(\vec{x};\,\boldsymbol{\mu}) &= z_N(\vec{x};\,\boldsymbol{\mu})+\eta_M(\vec{x};\,\boldsymbol{\mu}) \\
    &= \sum_{n=1}^N\alpha_n(\boldsymbol{\mu})\cdot \zeta_n(\vec{x})+\sum_{m=1}^M \theta_m(\boldsymbol{\mu})\cdot g_m(\vec{x})
\end{split}
    \label{eqn: PBDW-estimate}
\end{equation}

The update term $\eta_M$ accommodates non-modelled physics and unanticipated uncertainty, and it belongs to the functional space $U_M\subset \mathcal{U}$ spanned by $\{g_m\}_{m=1}^M$, hence this term is used to increase the knowledge of the prediction by taking into account what the model misses.

\subsubsection{Offline Phase: Construction of the Spaces}

The aim of the offline phase consists in the generation of a suitable basis for the background $Z_N$ and update $U_M$ space, namely the basis functions $\{\zeta_n\}_{n=1}^N$ and $\{g_m\}_{m=1}^M$ must be defined according to some optimality criteria. Different choices can be adopted for the selection of $\{\zeta_n\}_{n=1}^N$: the Proper Orthogonal Decomposition (POD) \cite{cordier_proper_2008, quarteroni2015reduced, brunton_data-driven_2022}, the Weak-Greedy \cite{prudhomme_mathematical_2002} or an hybrid version POD-greedy \cite{haasdonk_bernard_convergence_2013}. In this work, the first technique is used because is recognised as the state-of-the-art in dimensionality reduction. Given the training snapshots $\{u^{(i)}=u(\vec{x};\,\boldsymbol{\mu}_i\}_{i=1}^{N_s}\}$, the basis functions\footnote{In the context of POD, they may be referred to as spatial or POD modes.} depend on the eigenvalues-eigenvectors of the correlation matrix $\mathbb{C}\in\mathbb R^{N_s\times N_s}$, defined as
\begin{equation}
    \mathbb{C}_{ij} = \left(u^{(i)}, \,u^{(j)}\right)_{\mathcal{U}} \qquad i,j=1, \dots, N_s
\end{equation}
from which its eigenvalues $\{\lambda_k\}$ and eigenvectors $\{\eta_k\}$ are extracted. Then, the POD modes $\{\zeta_n\}_{n=1}^N$ can be computed using the following
\begin{equation}
    \zeta_n(\vec{x}) = \frac{1}{\sqrt{\lambda_n}}\,\sum_{i=1}^{N_s}\eta_{n,i}\cdot u^{(i)}(\vec{x})
\end{equation}

The POD modes defined in this way are ortho-normal with respect to the scalar product in $\mathcal{U}$ \cite{quarteroni2015reduced}: the most common choice is the $L^2$ scalar product.

The generation of a basis for the update space $U_M$ is directly linked to the search for the optimal positions of sensors, defined as linear functionals as in Section \ref{sec-geim}: in fact, the basis elements $\{g_m\}$ are taken as the Riesz representation\footnote{Let $v\in\mathcal{U}'$ be a linear functional onto the space $\mathcal{U}$, its Riesz representation $\mathcal{R}_\mathcal{U} v$ is defined as an element $g\in\mathcal{U}$ such that
\begin{equation}
    \left(\mathcal{R}_\mathcal{U}v, \varphi\right)_\mathcal{U} = 
    \left(g, \varphi\right)_\mathcal{U} = v(\varphi) \qquad \forall \varphi \in \mathcal{U} 
\end{equation}
given $(\cdot,\,\cdot)_\mathcal{U}$ as the inner-product in $\mathcal{U}$.} of the functionals $\{v_m\}$, thus it is sufficient to find a suitable procedure for the sensor positioning.

To find an optimal criterion for their placement, the \textit{a priori} PBDW error estimate must be used as driving choice \cite{maday_pbdw_2015, taddei_model_2016}: to ensure that the sensors are placed where the most important evolution takes place for the spatial modes $Z_N$ and that the update space represents the unknown compared to the \textit{background} model, a stability-approximation balancing algorithm (Algorithm \ref{Algo: SGREEDY}), called \textbf{SGreedy} is used \cite{maday_pbdw_2015, maday_adaptive_2019}.

\begin{algorithm}[ht]
\textbf{Input}\\
$\quad$ Reduced Space $Z_{N}$\;
$\quad$ Max Number of Sensors $M$\;
$\quad$ Library of sensors $\Upsilon$\;
\textbf{Output}\\
$\quad$ Update Space $U_M$\; $ $\\
\textbf{First Iteration} $m=1$\\
$\quad$ Search $v_1 = \argmax\limits_{v\in\Upsilon} v(\zeta_1)$\;
$\quad$ Set $g_1 = \mathcal{R}_\mathcal{U}v_1$ and $U_1 = \text{span}\{g_1\}$\;
\While {$m\leq M$}{
$\quad$ Set $n=\min\{N, m\}$\;
$\quad$ Compute $\beta_{n,m}=\inf\limits_{w\in Z_n}\sup\limits_{\varphi \in U_m} \frac{(w, \varphi)_\mathcal{U}}{\|w\|\cdot \|\varphi\|}$\;
$\quad$ Compute the least stable mode $w_{inf}=\arginf\limits_{w\in Z_n}\sup\limits_{\varphi \in U_m} \frac{(w, \varphi)_\mathcal{U}}{\|w\|\cdot \|\varphi\|}$\;
$\quad$ Compute the associated supremizer $ w_{sup,m} = \Pi_{U_m} w_{inf}$\;
$\quad $ Search $v_m = \argmax\limits_{v\in\Upsilon} v\left(w_{inf}-w_{sup,m}\right)$\;
$\quad$ Set $g_m = \mathcal{R}_\mathcal{U}v_m$ and $U_{m+1} = \text{span}\{U_m, g_m\}$\;
$\quad$ $m = m+1$ }
\caption{SGreedy for sensor selection}
\label{Algo: SGREEDY}
\end{algorithm}

The term $\beta_{N,M}$ is known as the \textit{inf-sup} constant and it is directly involved in the \textit{a priori} error estimate of the PBDW \cite{maday_parameterized-background_2014, maday_pbdw_2015, taddei_model_2016} (see Theorem \ref{theo: PBDWerror}). From a practical point of view, its calculation can be related to an eigenvalue problem \cite{maday_parameterized-background_2014}, and this equivalency is used within the SGreedy algorithm \ref{Algo: SGREEDY}.

\subsubsection{Online Estimation}

As for GEIM, the data are collected from the system \eqref{eqn: sec2-online-measure} to find the coefficients $\boldsymbol{\alpha}\in\mathbb{R}^N$ and $\boldsymbol{\theta}\in\mathbb{R}^M$ for the state estimation. They are the solution of the following linear system of dimension $(N+M)$ \cite{maday_yvon_pbdw_2015}
\begin{equation}
	\left[ 
	\begin{array}{cc}
		\mathbb{A}+\xi M\mathbb{I}& \mathbb{K}  \\ 
		\mathbb{K}^T & 0
	\end{array}
	\right] \cdot
	\left[ 
	\begin{array}{c}
		\boldsymbol{\alpha} \\ \boldsymbol{\theta}
	\end{array}
	\right]   =
	\left[ 
	\begin{array}{c}
		\vec{y} \\ \vec{0}
	\end{array}
	\right]\qquad
 \text{ given }\left\{
 \begin{array}{l}
      \mathbb{A}_{mm'}=\left(g_m,\,g_{m'}\right)_{\mathcal{U}} \\
      \mathbb{K}_{mn}=\left(g_m,\,\zeta_{n}\right)_{\mathcal{U}}
 \end{array}
 \right.
 \label{eqn: sec2-pbdw-online}
\end{equation}
with $m,m'=1,\dots M$ and $n=1,\dots, N$. The overall online procedure is summarised in Algorithm \ref{Algo: PBDW-Online}.

\begin{algorithm}[ht] 
	\textbf{Input}\\
	$\quad$ Basis functions $\{\zeta_{1}(\mathbf{x}),...,\zeta_{M}(\mathbf{x})\}$\;
	$\quad$ Basis sensors $\{v_{1},...,v_{M}\}$\;
	$\quad$ Riesz Representation of the sensors $\{g_{1},...,g_{M}\}$\;
	\textbf{Output}\\
	$\quad$ State Estimation $u^\star = z_N+\eta_M$\; $ $\\
	\textbf{Online Estimation}\\
	$\quad$ Acquisition of Experimental Data\\
	$\qquad \left\{y_m = v_m(u^{\text{true}})+\epsilon_m\right\}_{m=1}^M$\;
    $\quad$ Set the parameter $\xi$\;
    $\quad$ Solve for $\boldsymbol{\alpha}$ and $\boldsymbol{\theta}$ from the linear system \eqref{eqn: sec2-pbdw-online}\;
    $\quad$Estimate the state of the system\\
	$\qquad \displaystyle u^{\text{true}} \simeq u^\star = z_N+\eta_M = \sum_{n=1}^N\alpha_n\cdot \zeta_n+\sum_{m=1}^M \theta_m\cdot g_m$
	\caption{PBDW - Online Stage}
	\label{Algo: PBDW-Online}
\end{algorithm}

As for (TR-)GEIM, the offline and online phases of PBDW are implemented within the Python package \href{https://github.com/ERMETE-Lab/ROSE-pyforce}{pyforce}.

\subsection{Data-Driven Model Bias Correction}

The PBDW has been conceived to provide as estimated based on the \textit{projection-by-data}, introducing the update term $\eta_M$ to correct model bias \cite{haik_real-time_2023}, whereas the GEIM algorithm can introduce information through the interpolation condition and it is more suited when the number of sensors is limited \cite{riva_hybrid_2023-1}. Nevertheless, the latter can be seen as a special case of the latter \cite{taddei_model_2016} with perfect data ($\xi=0$), if the reduced space and the update one are generated with the GEIM greedy procedure and the number of measurements is equal to the dimension of the reduced space. 

Two important results from the \textit{a priori} error analysis of \cite{maday_adaptive_2019} will be discussed to show that the error of the PBDW estimate is controlled. At first, the case with perfect measures, without noise, will be presented, and then the results will be extended to the case of noisy data.

\begin{theo}[PBDW \textit{a priori} error estimates for perfect measurements]\label{theo: PBDWerror}
    Let $\{y_m = v_m\left(u^{\text{true}}\right)\}_{m=1}^M$, let $\beta_{N,M}>0$ and let $U_M$ satisfy the unsolvency condition
    \begin{equation*}
        \{v_m(\psi)=0\Longleftrightarrow \psi=0\}_{m=1}^M\qquad \psi\in U_M
    \end{equation*}
    The following estimates hold
    \begin{equation*}
        \begin{aligned}
            \norma{ u^{\text{true}}-u^{\text{opt}}_{\xi=0}}
            & \leq \frac{1}{\beta_{N,M}} \inf\limits_{z\in Z_N} \inf\limits_{\eta\in{U}_M\cap Z_M^\perp}\norma{u^{\text{true}}-z-\eta}\\
            \norma{ u^{\text{true}}-u^{\text{opt}}_{\xi=0}}
            & \leq \frac{\sqrt{4+6\norma{\mathcal{I}_M}^2_{\mathcal{L}(\mathcal{U})}}}{\beta_{N,M}} \inf\limits_{z\in Z_N} \inf\limits_{\eta\in{U}_M\cap Z_M^\perp}\norma{u^{\text{true}}-z-\eta}
        \end{aligned}
    \end{equation*}
    given $\mathcal{U}$ an Hilbert space, $\norma{\mathcal{I}_M}^2_{\mathcal{L}(\mathcal{U})}$ as the Lebesgue Constant in the $\mathcal{U}-$norm (it can be related to the inf-sup constant as in \cite{maday_generalized_2015,maday_adaptive_2019}) and $u^{\text{opt}}_{\xi=0}$ the solution to the PBDW statement with interpolation condition \cite{maday_parameterized-background_2014}.
\end{theo}

This result is referred to the case $\xi =0$ in which the reconstruction $u^\star$ interpolates the data, in particular as the \textit{inf-sup} constant increases (the update space $U_M$ becomes richer and richer) and the error decreases. However, it is important to retrieve a notion of convergence when $\xi >0$, i.e. for imperfect measures. By assuming the experimental data to be polluted by random noise (independent and distributed as a normal $\epsilon_m \sim \mathcal{N}(0,\sigma^2)$) and by defining the matrix $\mathbb{P}(\xi)$ as
\begin{equation}
    \mathbb{P} = \mathbb{P}(\xi) = \left[ 
	\begin{array}{cc}
		\mathbb{A}+\xi M\mathbb{I}& \mathbb{K}  \\ 
		\mathbb{K}^T & 0
	\end{array}
	\right]
\end{equation}

Let $\vec{u}^\star_\xi = [\boldsymbol{\alpha}_\xi, \boldsymbol{\theta}_\xi]^T$ and $\vec{u}^{\text{opt}}_{\xi=0}=[\boldsymbol{\alpha}^{\text{opt}}_{\xi=0}, \boldsymbol{\theta}^{\text{opt}}_{\xi=0}]^T$ be the coefficients of the linear expansion \eqref{eqn: PBDW-estimate} obtained from the solution of the algebraic problem \eqref{eqn: sec2-pbdw-online}, the following theorem hold:
\begin{theo}[PBDW \textit{a priori} error estimates for imperfect measurements]
    Supposing $\epsilon_m\sim\mathcal{N}(0,\sigma^2)$, the following hold
    \begin{equation*}
        \begin{aligned}
            \norma{\mathbb{E}[\vec{u}^\star_\xi] - \vec{u}^{\text{opt}}_{\xi=0}}_2&\leq \frac{\xi \,M}{s_{min}(\mathbb{P})}\norma{\boldsymbol{\theta}_{\xi=0}^{\text{opt}}}_2\\
            \mathbb{E}\left[\norma{\vec{u}^\star_\xi - \vec{u}^{\text{opt}}_{\xi=0}}_2^2\right]&\leq\left( \frac{\xi \,M}{s_{min}(\mathbb{P})}\norma{\boldsymbol{\theta}_{\xi=0}^{\text{opt}}}_2\right)^2+\sigma^2\cdot Tr\left(\mathbb{P}^{-1}\tilde{\mathbb{I}}\mathbb{P}^{-T}\right)
        \end{aligned}
    \end{equation*}
    given $\mathbb{E}[\cdot]$ the expected value operator,    $\tilde{\mathbb{I}} = \left[ \begin{array}{cc}
    \mathbb{I} & {0} \\ {0} & {0}
    \end{array}\right]$, $Tr(\cdot)$ indicating the trace of a matrix and $s_{min}(\cdot)$ the minimum singular value of a matrix.
\end{theo}

Thus, the error estimate is the sum of two contributions: the former involves the accuracy of the background space, measured and is monotonically increasing with $\xi$; the latter involves the accuracy of the measurements and is monotonically decreasing with $\xi$ according to \cite{maday_adaptive_2019}. This implies that there exists an optimal value of $\xi$ depending on the ratio between “measurement inaccuracy” and “model inaccuracy”. The more the model is accurate, the higher $\xi$ will be; on the other hand, the higher the accuracy of the measurements, the lower $\xi$ is. Since the selection of this parameter is a critical part of the online stage, proper adaptive strategies should be used exploiting cross-validation \cite{maday_yvon_pbdw_2015, maday_adaptive_2019}.

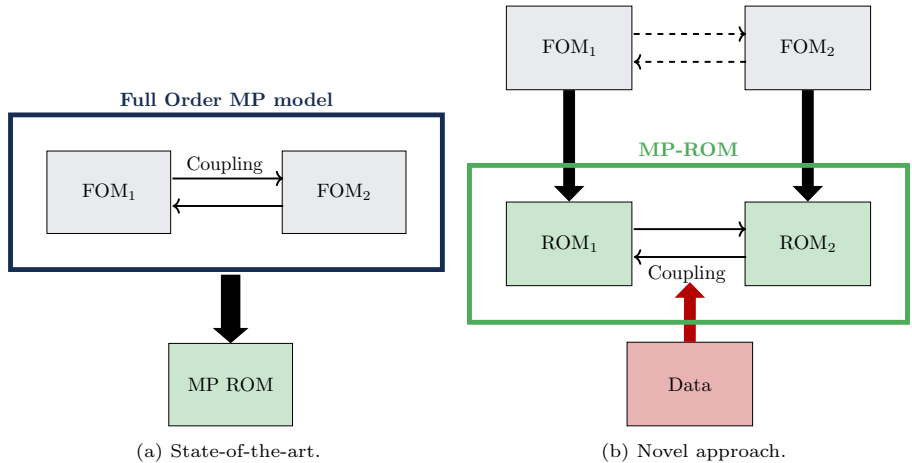
\begin{figure}[htbp]
    \centering
    \begin{subfigure}{0.495\textwidth}
		\centering
        \resizebox{1.\linewidth}{!}{
        \input{schemes/mp_rom_start_art}
        }
        \caption{State-of-the-art.}
        \label{fig: sota-mp-rom}
    \end{subfigure}\hfill
    \begin{subfigure}{0.495\textwidth}
        \centering
        \resizebox{1.\linewidth}{!}{\input{schemes/mp_rom_novel}
        }
        \caption{Novel approach.}
        \label{fig: novel-rmp}
    \end{subfigure}
    \caption{Reduced Order Modelling for Multi-Physics systems: state-of-the-art and proposed approach.}
    \label{fig: rom-mp}
\end{figure}

The possibility of updating/correcting the information coming from the mathematical model is very interesting in the nuclear field, especially when MP modelling is involved. In the last few years, nuclear reactors modelling is being shifted from very accurate single physics codes to MP framework by gathering all the governing equations inside the same environment \cite{fiorina_gen-foam_2015, demaziere_6_2020}. The latter possibility has become possible due to increasing computational capabilities\footnote{Even though, powerful machine are still required.}, however, in this way, all the previous validation of single-physics code is getting lost. In order to keep this aspect and the experience gained on single-physics code over the year, a novel approach can be proposed. Figure \ref{fig: sota-mp-rom} shows the state-of-the-art of Reduced Order Methods in the MP context considering two different physics involved: a fully coupled FOM is trained to generate the snapshots and then a representation in the reduced coordinate system is defined for quick and efficient evaluations \cite{lassila_model_2014}. On the other hand, Figure \ref{fig: novel-rmp} shows the approach proposed in this paper: assuming to have at disposal two single physics codes which are not able to directly communicate, some training snapshots are generated introducing the coupling in a weaker/approximated way, then the reduced representation is derived from each set of snapshots so that during the online phase it is enriched by the information of the data, which is intrinsically MP, sharing additional information about the coupling.

It is crucial to enforce some information of the coupling to the training snapshots, to have a proper reduced representation which is aware of the coupling effects and hence can incorporate the coupling information of the data\footnote{If the FOMs were completely de-coupled, the ROM cannot be used to extrapolate the coupling from the data unless a lot of data are used, namely a lot of measurements with $\vec{y}\in\mathbb{R}^M$ of high dimension. This is unfeasible because the number of sensors is typically limited \cite{riva_hybrid_2023}.}. 

The possibility of correcting the prediction of the offline phase will be investigated in two contexts:
\begin{enumerate}
    \item the single-physics codes cannot be directly communicating and the snapshots are generated in an approximated way using a weak Picard iteration (Section \ref{sec: aFOM-coupling});
    \item the single-physics code can directly communicate, however, the coupling function may not be known and hence, it is approximated with a linear function (Section \ref{sec: lFOM-coupling}).
\end{enumerate}
In the end, the FOM considered as the \textit{truth} from which data are generated will be described in Section \ref{sec: FOMcoupling}.

\section{Multi-Physics Case Studies}\label{sec: MP-case-study}

In nuclear reactors, the most important coupling is the one between Neutronics and Thermal-Hydraulics (N-TH) \cite{demaziere_6_2020} as they influence each other through feedback effects: the power generated by the fission event induces a temperature profile, which changes the material properties in terms of densities and microscopic cross sections \cite{DuderstadtHamilton}. The modelling of this temperature feedback is typically done through correlations as in \cite{Aufiero2014}, by learning the map through Machine Learning techniques \cite{PHYSOR24_MC_ML} or with a more advanced coupling exploiting nuclear libraries as in \cite{CastagnaIntroini_2022}. 

In this work, two test cases with coupled neutronics-thermal hydraulics are analysed: the IAEA PWR 2D \cite{argonne_book, Milonga_code} and the TWIGL2D reactor \cite{twigl-Hageman, twigl-Yasinsky}. These are simple problems to be implemented in a numerical solver, however, they include the main features of the N-TH coupling in commercial nuclear reactors and they are used to assess the capabilities of new neutronics solvers. In the following section, the governing equations will be presented, followed by the coupling strategies to generate the \textit{truth} and the training snapshots.

\subsection{Governing Equations}

The neutronics is modelled using the multi-group neutron diffusion equations \cite{DuderstadtHamilton, hebert_applied_2009}, a commonly used approximation to model the transport of neutrons inside multiplying media. Let $\phi_g$ be the neutron flux of the $g-$th group and $c_j$ be the $j-$th precursors groups, satisfying the following system of PDEs ($g=1, \dots, G$ and $j=1, \dots, J$)
\begin{equation}
    \left\{
        \begin{aligned}
            \frac{1}{v_g}\dpart{\phi_g}{t}=&+\nabla\cdot(D_g\nabla \phi_g)-\left(\Sigma_{a,g}+\sum_{g'\neq g}\Sigma_{s,g\rightarrow g'} + D_g B_{z,g}^2\right)\phi_g\\
            &+\sum_{g'\neq g}\Sigma_{s,g'\rightarrow g}\phi_{g'}+\chi_g\left(\frac{1-\beta}{k_{\text{eff}}}\sum_{g'}\nu_{g'}\Sigma_{f,g'}\phi_{g'}+\sum_{j=1}^J\lambda_jc_j\right)\\
            \dpart{c_j}{t} =&  \frac{\beta_j}{k_{\text{eff}}}\sum_{g}\nu_g\Sigma_{f,g}\phi_{g} - \lambda_jc_j
        \end{aligned}
    \right.
    \label{eqn: mg-neutron-diffusion}
\end{equation}
completed with vacuum and symmetry boundary conditions for the fluxes, with the main quantities listed in Table \ref{tab: Neutronic_quantities}. 

\begin{table}[htbp]
	\centering
	\begin{tabular}{|c|l|}
        \toprule
        \textbf{Symbol} & \textbf{Meaning}\\\midrule
		$\phi_g$ & Group flux\\
		$v_g$  & Group neutron velocity\\
		$D_g$ & Diffusion coefficient of group $g$\\
		$\Sigma_{a,g}$ & Macroscopic absorption cross section of group $g$\\
		$\Sigma_{s,g\rightarrow g'}$ & Macroscopic scattering cross section from $g$ to $g'$\\
		$B_{z,g}$ & Axial buckling of group $g$\\
		$\chi_{g}$ & Fission spectrum of group $g$\\
		$k_{\text{eff}}$ & Effective multiplication factor\\
		$\beta$ & Total delayed neutron fraction \\
		$\Sigma_{f,g}$ & Macroscopic fission cross-section of group $g$\\
		$\nu_g$& Average neutrons produced per fission in group $g$ \\
		$c_j$& Delayed neutron precursor concentration of the $j^{\text{th}}$ group\\
		$\beta_j$& Delayed neutron fraction of the $j^{\text{th}}$ precursor group\\
		$\lambda_j$& Decay constant of the $j^{\text{th}}$ precursor group\\
		\bottomrule
	\end{tabular}
	\caption{Overview of quantities involved and terms in the multi-group neutron diffusion approximation.}
	\label{tab: Neutronic_quantities}
\end{table}

On the other hand, the thermal part involves only the heat diffusion without fluid dynamics\footnote{This is done to have a simpler model onto which the pros and cons of the model correction for MP problems can be easily observed. In the future, more complex models will be considered.}: the temperature $T$ is governed by the unsteady heat equation
\begin{equation}
    \rho c_p \dpart{T}{t} = \nabla\cdot (k\nabla T)+q'''
    \label{eqn: heat-equation}
\end{equation}
in which $\rho$ is the density, $c_p$ the specific heat capacity at constant pressure, $k$ the thermal conductivity and $q'''$ the power density coming from the fission event. The PDE is completed with fixed temperature (Dirichlet BC) or convective heat transfer (Robin BC) and symmetry boundary conditions. These equations are coupled through the temperature dependence of the absorption cross sections and the diffusion coefficients, through the fission power density
\begin{equation}
    \left\{
        \begin{aligned}
            D_g &= f_D\left(T;\,T^{\text{ref}}, D_g^{\text{ref}},\boldsymbol{\gamma}_D\right) & g= 1, \dots G\\
            \Sigma_{a,g} &= f_a\left(T;\,T^{\text{ref}}, \Sigma_{a,g}^{\text{ref}},\boldsymbol{\gamma}_a\right) & g= 1, \dots G\\
            q''' &= P_0\cdot \sum_{g=1}^G \Sigma_{f,g}\phi_g
        \end{aligned}
    \right.
    \label{eqn: coupling}
\end{equation}
given $P_0$ a normalisation constant and $f_D$ and $f_a$ suitable functions of the temperature dependent on some coupling parameters $\boldsymbol{\gamma}_D$ and $\boldsymbol{\gamma}_a$.

The governing equations have been discretised in space using the Finite Element Method \cite{quarteroni_numerical_2016}, exploiting the FEniCSx library for Python, whereas implicit Euler schemes are used to advance in time; for further details about the theory of finite elements applied to neutron transport refer to \cite{ackroyd_finite_1997} and for this specific problem see \ref{app: FE-MG-diffusion}. In the following subsections, three different ways of coupling the equations will be presented.

\subsubsection{High-Fidelity Full Order Model}\label{sec: FOMcoupling}

The solution to the fully coupled system of PDE represents the \textit{truth}, i.e. $u^{\text{true}}$, from which synthetic data are generated for the online phase. The multi-group diffusion equations \eqref{eqn: mg-neutron-diffusion} and the heat equation \eqref{eqn: heat-equation} are non-linear due to the dependence of the absorption cross sections and diffusion coefficients on the temperature field through $f_a$ and $f_a$ respectively. A segregated approach has been selected in this work \cite{demaziere_6_2020}, hence at each time step (recall that implicit Euler is used) a non-linear system of equations should be solved. In order to retrieve linearity at each time step, it has been chosen to adopt a semi-implicit treatment of the non-linear terms \cite[Chapter~16]{quarteroni_numerical_2016}: let the superscript $^{(n)}$ be used to indicate the evaluation of a function at time $t=t^{(n)}$, the non linear terms are written as follows at $t^{(n+1)}$:
\begin{equation}
\left\{
    \begin{aligned}
        \nabla\cdot \left(D_g^{(n)}\nabla \phi_g^{(n+1)}\right) \qquad & \text{ with } D_g^{(n)} = f_D\left(T^{(n)};\,T^{\text{ref}}, D_g^{\text{ref}},\boldsymbol{\gamma}_D\right)  \\
        \Sigma_{a,g}^{(n)} \phi_g^{(n+1)} \qquad & \text{ with } \Sigma_{a,g}^{(n)} = f_a\left(T^{(n)};\,T^{\text{ref}}, \Sigma_{a,g}^{\text{ref}},\boldsymbol{\gamma}_{a,g}\right) 
    \end{aligned}
\right.
\end{equation}
This approximation scheme assures a first-order convergence with the time step size \cite{quarteroni_numerical_2016}.

\subsubsection{Approximated Full Order Model: External Picard Coupling}\label{sec: aFOM-coupling}

The single-physics codes for the approximated FOM are the neutronic (N) and thermal (TH) solver, with the assumption that they cannot interact while advancing in time, meaning that the output of one cannot be directly used as input to the other. However, it is necessary to enforce some information about the coupling effects inside the system to aid the correction at the reduced level, as explained in Figure \ref{fig: novel-rmp}; hence, the temperature profile and the power density part must be given as "boundary conditions" to the neutronics and the heat equations respectively. Nevertheless, assuming that both of them can read the output of the \textbf{pyforce} library, the POD with Interpolation (POD-I) \cite{PODi_demo} can be exploited to create the "boundary condition" per each solver: this can be considered as a surrogate modelling for external physics. In this way, for instance, the power density dependent on time is represented in a reduced coordinate system and it can be inserted as external input without increasing the computational burden of the solution\footnote{The spatial behaviour is caught by the POD modes $\{\zeta_n\}$, whereas the time/parametric dependence is hidden in the modal coefficients $\{\alpha_n\}$. The POD-I algorithm requires a map between $\boldsymbol{\mu}$ and the coefficients $\{\alpha_n\}$, using interpolants as in \cite{riva_hybrid_2023-1} or machine learning techniques as in \cite{PODI-GP}. The former strategy is adopted in this work}.

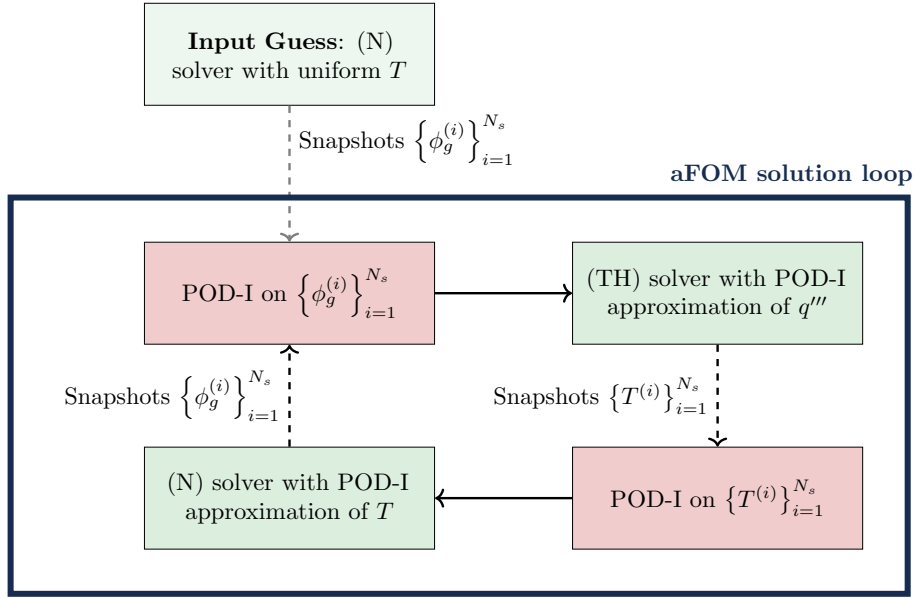
\begin{figure}[htbp]
    \centering
    \resizebox{1.\linewidth}{!}{\input{schemes/pod-training}}
    \caption{Approximated Full Order Model solution scheme.}
    \label{fig: POD-training-aFOM}
\end{figure}

The idea of the coupling enforcement explained above is reported in Figure \ref{fig: POD-training-aFOM}: in this way, the solvers do not interact directly but only through the reduced representation. 
In this work, only 2 only 2 iterations of the weak coupling loop are performed; this scheme is used to obtain a low-fidelity solution acting as an initial guess which will be corrected when the knowledge of the data is introduced in the online phase. This was of proceeding allows the save of computational resources, which would be required to get a fully converged solution by a Picard iteration \cite{demaziere_6_2020}. This approach of generating the snapshots will be referred to as aFOM.

This strategy has been chosen due to its similarity with the Picard iteration: the POD-I is crucial to generate a surrogate model serving as a boundary condition to the neutronic and/or thermal-hydraulic solver. Moreover, this way of proceeding keeps the solvers separated without direct communication to mimic what is typically encountered in the nuclear world. In the end, it should be pointed out that this strategy can be computationally expensive since the number of iterations may become high; in the future, the "boundary conditions" for each solver can be imposed using surrogate modelling techniques coming directly from data.

\subsubsection{Linearised Full Order Model: Linear Coupling Dependence}\label{sec: lFOM-coupling}

Assuming now that the single-physics codes can be directly interfaced as in the FOM (Section \ref{sec: FOMcoupling}), this second way of generating the snapshots considers an approximated functional dependence for the absorption cross sections and diffusion coefficients through a linear approximation
\begin{equation}
    \left\{
        \begin{aligned}
            f_D\left(T;\,T^{\text{ref}}, D_g^{\text{ref}},\boldsymbol{\gamma}_D\right) &\approx m_{D_g}\cdot T + q_{D_g} \\
            f_a\left(T;\,T^{\text{ref}}, \Sigma_{a,g}^{\text{ref}},\boldsymbol{\gamma}_a\right) &\approx m_{\Sigma_{a,g}}\cdot T+ q_{\Sigma_{a,g}}\\
        \end{aligned}
    \right.
    \label{eqn: lin-coupling-fd-fa}
\end{equation}
given $m$ as the slope of the line and $q$ its intercept. This model will be labelled as LcFOM standing for Linearised Coupling Full Order Model. This approximation is introduced to analyse whether or not an MP model with this coupling can be corrected using data coming from the \textit{ground-truth}, i.e. the solution of the fully-coupled MP model.

\subsection{IAEA 2D PWR benchmark}\label{sec-IAEA2D-case-study}

\begin{figure}[htbp]
    \centering
    \includegraphics[width=0.75\linewidth]{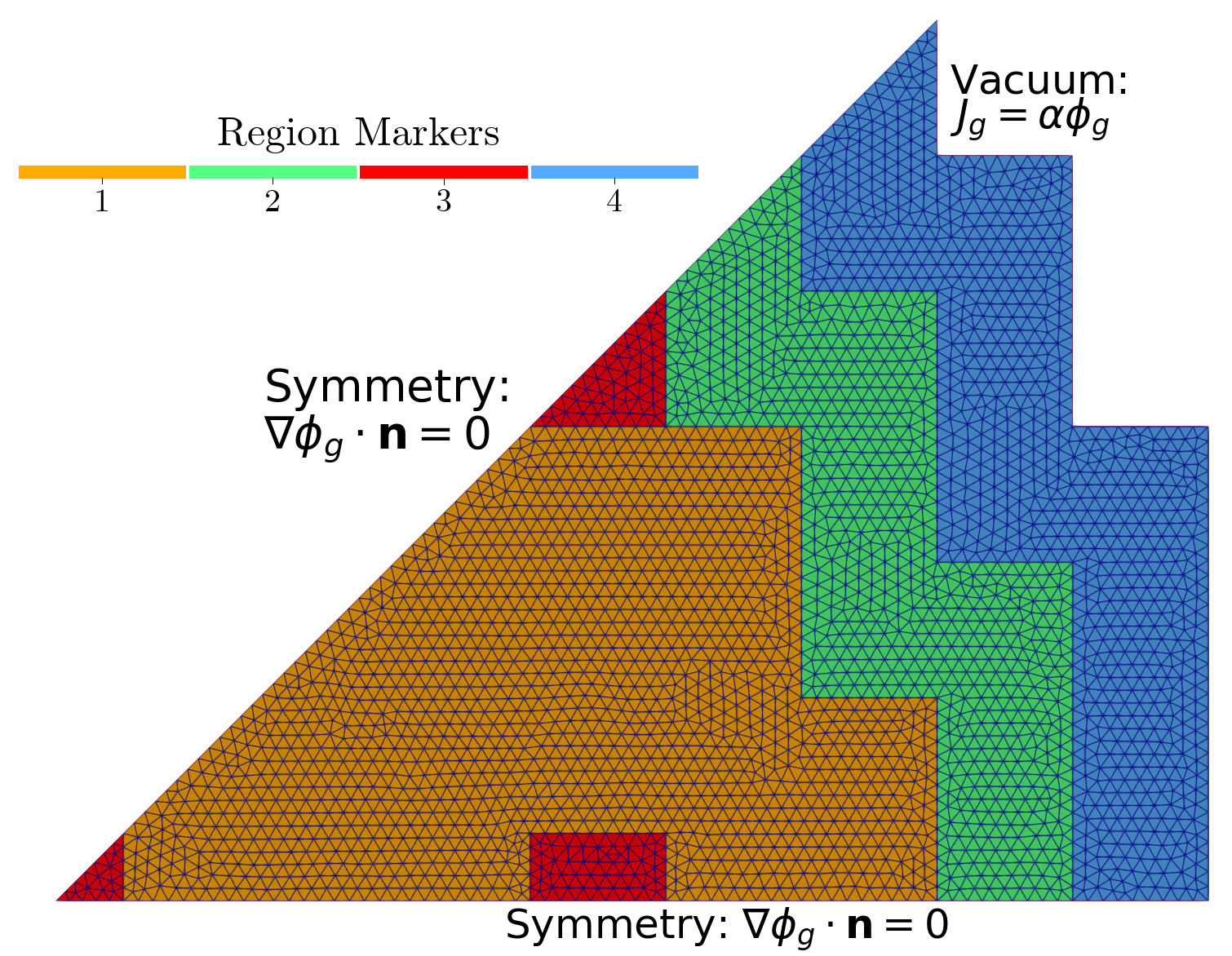}
    \caption{Regions, boundary conditions and mesh of the IAEA 2D PWR reactor.}
    \label{fig: pwr2d-regions}
\end{figure}

This case study is based on the ANL11-A2 benchmark from \cite{argonne_book}, which has been extended to time-dependent problems with thermal feedback to simulate an MP system. The geometry is reported in Figure \ref{fig: pwr2d-regions} and it consists of a one-octave of a PWR core in which there are represented 4 regions with different materials: two different fuels in regions 1 and 2, a homogeneous mixture of fuel 2 and control rod in region 3 and the reflector region in 4. 

For this case, there are $G=2$ energy groups\footnote{In the paper, the following convention is used: $g=1$ is referred to the fast flux, whereas $g=2$ is the thermal one.} and $J=6$ precursor groups; the transient to be simulated is a step insertion of reactivity obtained by reducing the absorption cross sections of both groups in region 3 as
\begin{equation}
    \Sigma_{a,g}^{\text{ref}}(\vec{x}, t) = \Sigma_{a,g}^{0}(\vec{x})\cdot 0.9\, 
    {H}(t)\qquad\qquad g = 1, \dots , 2
\end{equation}
given ${H}(t)$ the Heaviside step function. The following notation for the superscript is adopted: $^{\text{ref}}$ indicates the value at any time $t$ with the temperature equal to the reference value (no feedback); whereas $^0$ is used to indicate the value at the initial time with the temperature equal to the reference one. In the end, for this problem the following dependencies of the material properties on the temperature $T$ are considered ($g={1,2}$)
\begin{equation}
\begin{split}
D_g(\vec{x},t) &= D_g^{\text{ref}}(\vec{x}) + \gamma_{D,g}\cdot \ln \frac{T(\vec{x},t)}{T^{\text{ref}}}\\
\Sigma_{a,g}(\vec{x},t)& = \Sigma_{a,g}^{\text{ref}}(\vec{x},t) + \gamma_{a,g}\cdot \ln \frac{T(\vec{x},t)}{T^{\text{ref}}}
\end{split}
\end{equation}
given $T^{\text{ref}}=600\,$K, $\gamma_{D,1}=\gamma_{D,2} = 3\cdot 10^{-3}\,$cm and $\gamma_{a,1} = \gamma_{a,2} = 2.034\cdot 10^{-3}\,$cm$^{-1}$. All the other physical parameters of this problem are reported in \ref{app: pwr2d-params}.

\subsection{TWIGL2D benchmark}\label{sec-TWIGL2D-case-study}

\begin{figure}[htbp]
    \centering
    \includegraphics[width=0.75\linewidth]{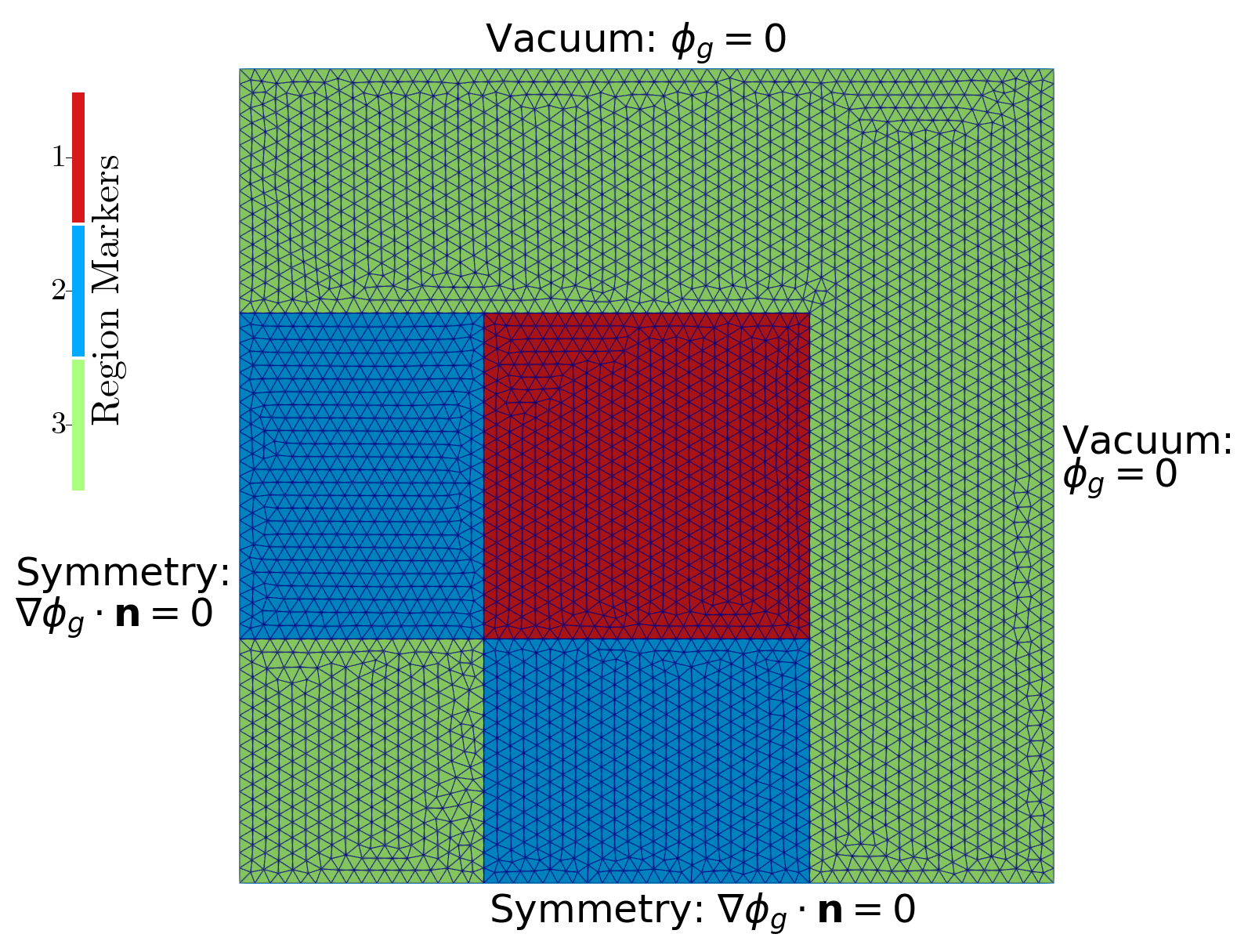}
    \caption{Regions of the TWIGL2D reactor.}
    \label{fig: twigl-regions}
\end{figure}

The two-dimensional TWIGL reactor kinetics problem was first presented in \cite{twigl-Hageman, twigl-Yasinsky} and it is a useful benchmark to assess the reliability and the efficiency of numerical solvers for transient neutronics \cite{Zachary2016, Zhang2014}. The system is composed of two different fuel regions (Figure \ref{fig: twigl-regions}): the seed regions 1 and 2 and the blanket zone labelled as 3. The former regions share the same material properties, however, they differ in the reactivity inserted during the transient. For this problem, $G=2$ energy groups are considered (1 = fast, 2 = thermal) and $J=1$ group of precursors.

The original benchmark problem involved only kinetics, hence no thermal feedback was considered: a positive ramp and a positive step in reactivity were given to the system by reducing the absorption cross-section of group 2 in region 1\footnote{The FEM solver, implemented in Python using FEniCSx, used in the MP problem has been previously validated against benchmark data. This analysis is not reported as it is out of the scope of the work, this analysis is not reported.}; in this paper, the following transient is studied\footnote{The superscript $0$ indicates that the temperature is equal to the reference value and that $t=0$, i.e. the beginning of the transient.} for $t\geq 0$
\begin{equation}
    \Sigma^{\text{ref}}_{a,2}(\vec{x},t) = \Sigma_{a,2}^{0}(\vec{x})\cdot \left[ 
    {H}(0.2-t)\cdot (1-0.11667\,t) + {H}(t-0.2)\cdot 0.97666 \right]
\end{equation}
for $\vec{x}\in\Omega_1$, given $\Omega_1$ as region 1 in Figure \ref{fig: twigl-regions} and ${H}(t)$ as the Heaviside step function. 

In the end, for what concerns the coupling functions two different cases are studied for the aFOM and LcFOM: the former considers that only $\Sigma_{a,1}$ can vary with the temperature with the following law (this case will be labelled as TWIGL2D-A)
\begin{equation}
    \Sigma_{a,1}(\vec{x},t) = \Sigma_{a,1}^{\text{ref}}(\vec{x})\cdot \left[1+\gamma_{a,1}\,\left(\sqrt{T(\vec{x},t)}-\sqrt{T^{\text{ref}}}\right)\right]
    \label{eqn: TWIGL-A-feedback}
\end{equation}
given $\gamma_{a,1}\in[10^{-3}, 10^{-2}]\;\left(\text{K}^{-1/2}\right)$ be one of the parameters in $\boldsymbol{\mu}$, in addition to time $t$, and $T^{\text{ref}} = 600\,$K; whereas, the latter involves the linearisation of the coupling function $f_D$ and $f_a$ defined as
\begin{equation}
\begin{split}
    f_{D} &= D^{\text{ref}}_g\cdot \left[1+\gamma_{D_g}\cdot \sin\left(2.75\cdot \frac{T-T^{\text{ref}}}{T^{\text{ref}}}\right)\right]\\
    f_{a} &= \Sigma_{a,g}^{\text{ref}}\cdot \left[1+\gamma_{\Sigma_{a,g}}\cdot \tanh\left(2\cdot \frac{T-T^{\text{ref}}}{T^{\text{ref}}}\right) \right]
\end{split}
    \label{eqn: TWIGL-B-feedback}
\end{equation}
given $T^{\text{ref}}=600\,$K, $\gamma_{D_1} = \gamma_{D_2} = 2\cdot 10^{-2}$ and $\gamma_{\Sigma_{a,1}} = \gamma_{\Sigma_{a,2}} = 4\cdot 10^{-2}$. The linearisation for them searches for the best linear fit in the temperature range $[600, 1200]$ using the \textit{scikit-learn} package \cite{scikit-learn}. This case will be labelled as TWIGL2D-B.
All the physical parameters of this problem are reported in \ref{app: twigl-params}.

\section{Numerical Results}\label{sec: num-res}

This section discusses the main numerical results of this work regarding the application of the TR-GEIM and PBDW to MP case studies presented in the previous section. These model correction algorithms combine the knowledge of mathematical modelling and evaluations on the system (i.e., data or measures) and the improvement depends strongly on the available data and on how they are collected: in particular, the uncertainty in the measurement procedure can be quite important especially in hostile environments as the one of nuclear reactors \cite{cammi_indirect_2023}. This work considers only synthetic data\footnote{This expression is used to identify a measurement vector coming from snapshots of a \textit{high-fidelity} simulation.} polluted by random noise, assumed to be distributed as a normal with null mean and standard deviation $\sigma = \left\{0.5\,\text{K} 0.01\,\frac{\text{n}}{cm^2\,s}, 0.01\,\frac{\text{n}}{cm^2\,s}\right\}$ for temperature $T$, fast flux $\phi_1$ and thermal flux $\phi_2$ respectively. Before entering the analysis of the main results, the available libraries of sensors must be specified: the general shape of the sensor has been given in Equation \eqref{eqn: sens-def}, the kernel function $\mathcal{K}\left(\|\vec{x}-\vec{x}_m\|_2, s\right)$ is a normalised Gaussian
\begin{equation}
    \mathcal{K}\left(\|\vec{x}-\vec{x}_m\|_2, s\right) =
    \frac{e^{\frac{-\|\vec{x}-\vec{x}_m\|_2^2}{2s^2}}}{\displaystyle\int_\Omega e^{\frac{-\|\vec{x}-\vec{x}_m\|_2^2}{2s^2}}\,d\Omega}
\end{equation}
in which the available positions $\vec{x}_m$ are the nodes of the numerical mesh\footnote{If the numerical mesh is composed of a lot of elements, the computational costs related to the selection of the positions of the sensors can become very high \cite{riva_hybrid_2023}. Thus, in this work it has been chosen to sample every 5 cells in the numerical mesh: this choice is coherent with the real world some positions can be unavailable if they are too close to other sensors.} and $s$ is the point spread taken as 1, given also the following property $\norma{\mathcal{K}}_{L^1(\Omega)}=1$. These functions are common for all the fields $(T, \phi_1,\phi_2)$.

As stated in Section \ref{sec-ROM}, ROM techniques are characterised by an offline and online phase. The former is devoted assessment of the reducibility of the problem \cite{lassila_model_2014, quarteroni2015reduced}, to the generation of the reduced space by defining the basis functions and to the generation of the basis sensors and hence to the search for their optimal location to maximise the amount of information extracted from the system. The latter phase consists of a quick state estimation, thus in the model correction phase; since the focus of this work is on the correction of MP models using data during the online phase, only the results of the latter phase will be presented. 

\subsection{Model Correction from Weakly Coupled FOM}

In this section, the numerical results related to the case in which the codes cannot directly communicate with one another is analysed, namely the model bias correction of aFOM using data coming from FOM. Before entering into the details, some definitions will be provided which will be used in this section. Let $\mathcal{P}_M\left[u(\cdot;\boldsymbol{\mu})\right]$ the reconstruction operator, either for (TR-)GEIM or PBDW, with $M$ measurements for the generic field $u(\cdot;\boldsymbol{\mu})$ and let $r\left[u(\cdot;\boldsymbol{\mu})\right](\vec{x}) = \left| u(\cdot;\boldsymbol{\mu}) - \mathcal{P}_M\left[u(\cdot;\boldsymbol{\mu})\right]\right|$ be the residual field. Moreover, the average absolute $E_M$ and relative $\varepsilon_M$ error measured in $\norma{\cdot}_{L^2(\Omega)}$ are defined as
\begin{equation}
\begin{split}
    E_M &= \frac{1}{\text{dim}(\Xi^{\text{predict}})}\sum_{\boldsymbol{\mu}\in\Xi^{\text{predict}}}\norma{r(\vec{x};\boldsymbol{\mu})}_{L^2(\Omega)}
\\
    \varepsilon_M &= \frac{1}{\text{dim}(\Xi^{\text{predict}})}\sum_{\boldsymbol{\mu}\in\Xi^{\text{predict}}}\frac{\norma{r(\vec{x};\boldsymbol{\mu})}}{\norma{u(\vec{x};\boldsymbol{\mu})}
          }_{L^2(\Omega)}
\end{split}
\label{eqn: error-defs}
\end{equation}
given $\Xi^{\text{predict}}\subset\mathcal{D}$ a subset of the parameter space with unseen data.

\subsubsection{IAEA 2D PWR Benchmark}

The problem described in Section \ref{sec-IAEA2D-case-study} considers the time $t$ as the only parameter $\boldsymbol{\mu}$, sampled in the parameter space $\mathcal{D}=[0, 2]\subset \mathbb{R}^+$ from which the training and predict/test set is defined $\Xi^{\text{train}} \subset [0, 1]$ and $\Xi^{\text{predict}} \subset (1, 2]$; in which the number of training snapshots used to generate the basis is $N_s = \text{dim}\left(\Xi^{\text{train}}\right)=101$, sampling the time $t$ every 0.01 seconds. These are used to generate the magic functions and sensors from the GEIM algorithm, and then the reduced space $Z_N$ for the PBDW is created with the POD: it has been decided to use $N=5$ because it provides a good approximation of the background model and assuring stability of the algorithm by keeping it low \cite{taddei_model_2016}. Moreover, the hyper-parameter $\xi$ has been tuned using a cross-validation procedure on a small validation set sampled in $\Xi^{\text{predict}}$: in particular, the reconstruction error $\varepsilon_M$ is optimised.
\begin{figure}[htbp]
    \centering
    \includegraphics[width=0.98\linewidth]{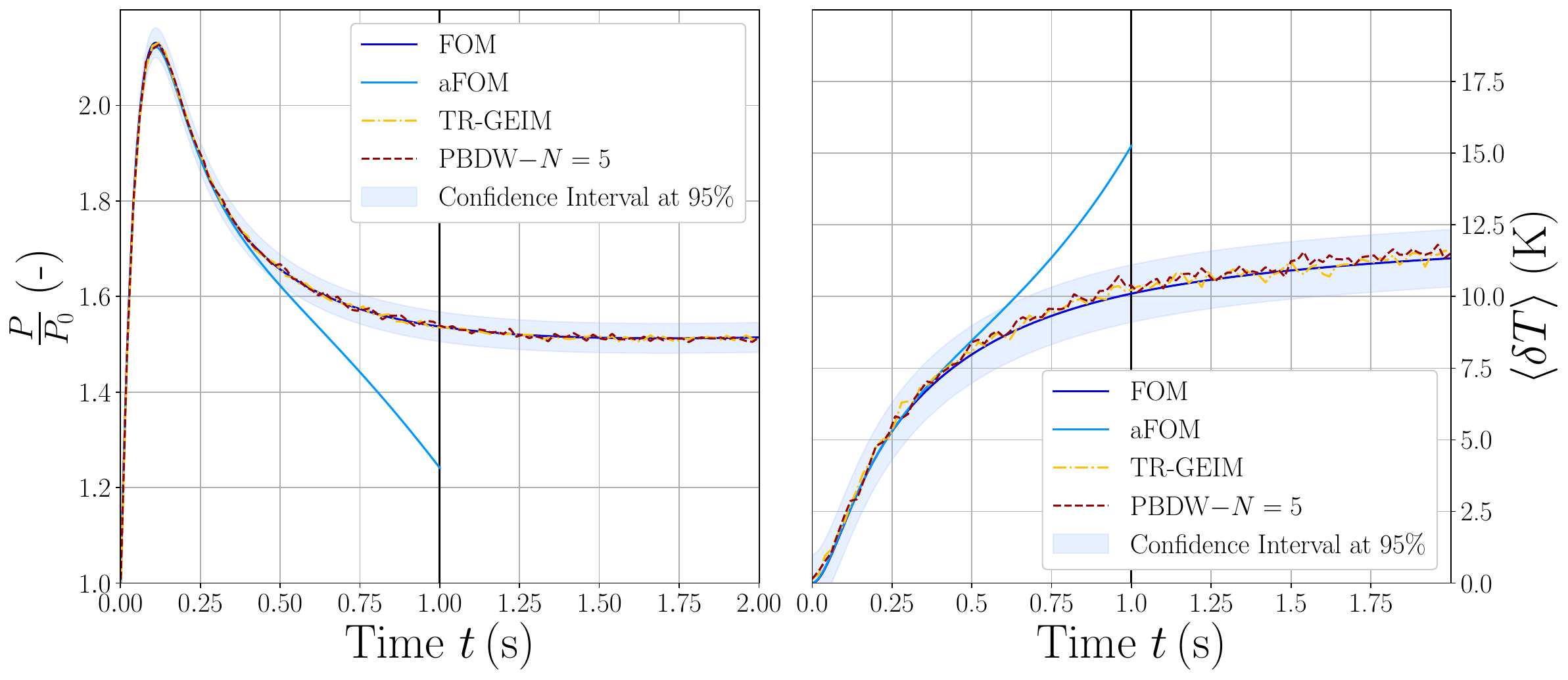}
    \caption{Line plots of the power (normalised with respect to the initial steady state) and the average temperature difference with respect to initial value for the IAEA 2D PWR case.}
    \label{fig: aFOM-ANL-LinePlot}
\end{figure}

Firstly, a comparison in terms of global output of interests is presented, with $M_{max}=15$ measurements, considering the overall power of the reactor and the average (in space) temperature difference compared to the initial state is reported in Figure \ref{fig: aFOM-ANL-LinePlot}, whose definitions are
\begin{equation}
    \begin{split}
        P(t) &= P_0\cdot \sum_{g=1}^G \int_\Omega \Sigma_{f,g}\phi_g(\vec{x},t)\,d\Omega\\
        \langle \delta T\rangle(t) &= \frac{\norma{T(\vec{x},t)}_{L^1(\Omega)} - \norma{T(\vec{x},0)}_{L^1(\Omega)} }{|\Omega|} 
    \end{split}
    \label{eqn: global-var-def}
\end{equation}
given $|\Omega| = \int_\Omega d\Omega$ as the measure of the spatial domain.

The output of the aFOM is much different with respect to the true value of the FOM, nevertheless, the TR-GEIM and the PBDW are able to retrieve the correct time evolution within the confidence interval\footnote{The error bands have been assessed using the Monte-Carlo method for Uncertainty Quantification.} at 95\%. 

\begin{figure}[htbp]
    \centering
    \includegraphics[width=0.98\linewidth]{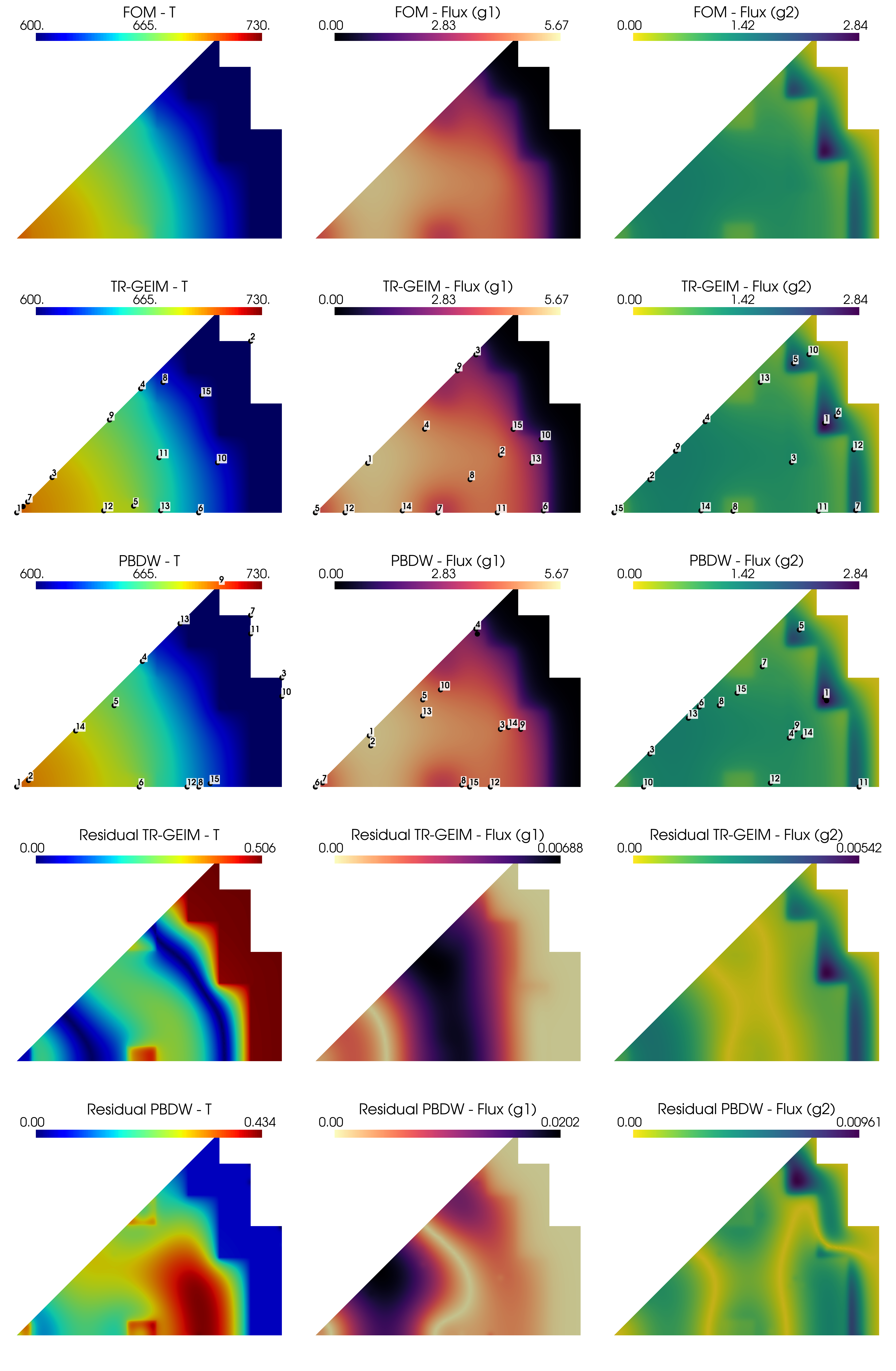}
    \caption{Contour plots of the temperature, fast flux and thermal flux fields (from left to right) for FOM, TR-GEIM and PBDW at final time $t=2$ seconds, and the residual field for the IAEA 2D PWR case.}
    \label{fig: aFOM-ANL-Contour}
\end{figure}

In addition to the comparison in terms of global output of interests, it is also important to assess what happens locally. Figure \ref{fig: aFOM-ANL-Contour} shows the temperature, fast flux (g1) and thermal flux (g2) at final time $t=2$ seconds, with $M_{max}=15$ measurements, considering the reconstruction through TR-GEIM and PBDW with $M=15$ sensors. Noisy data and inaccurate training snapshots do not limit the inference capabilities of these methods, on the contrary, they are able to correct the information provided in the offline training phase. The residual fields in the last two rows have low magnitudes highlighting the accuracy of the reconstructions.

\begin{figure}[htbp]
    \centering
    \includegraphics[width=1.\linewidth]{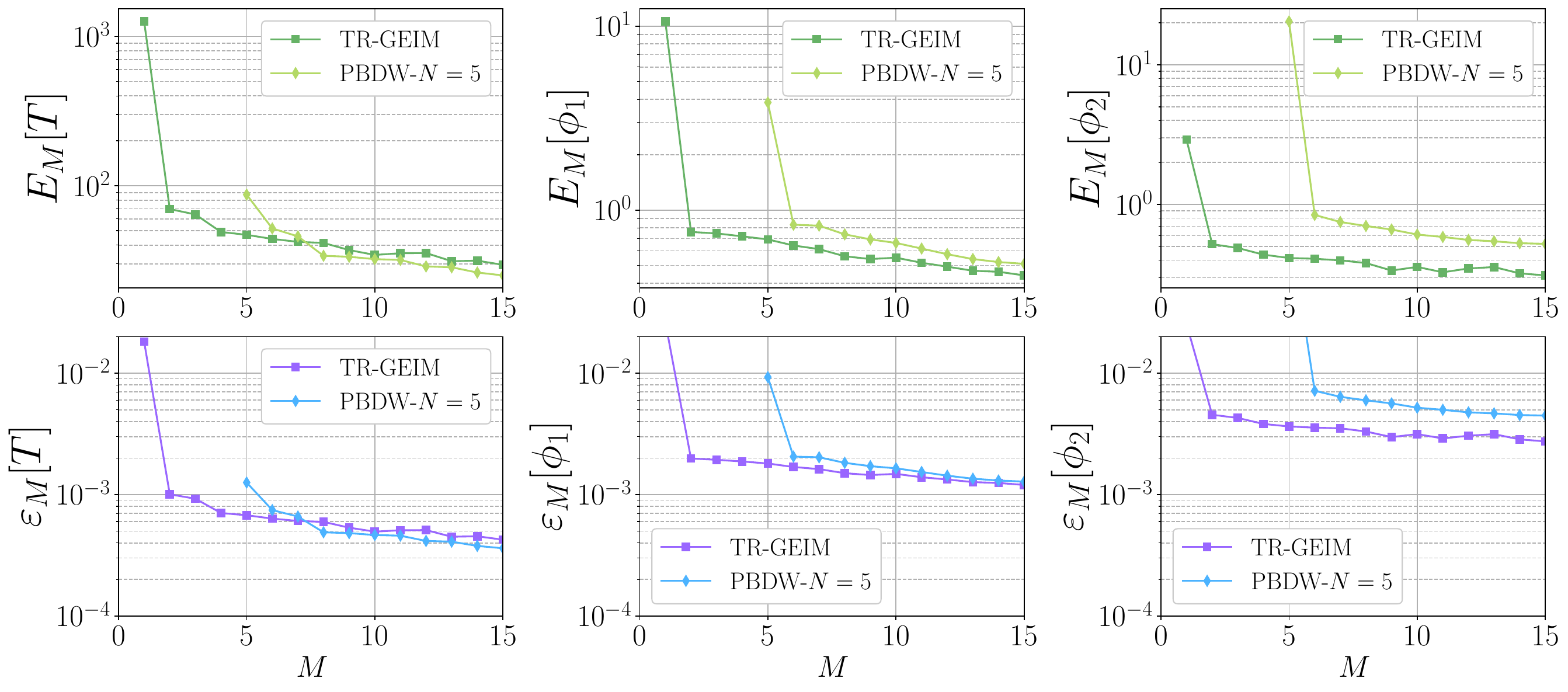}
    \caption{Average absolute and relative error measured in $\norma{\cdot}_{L^2(\Omega)}$, as in Eq. \eqref{eqn: error-defs}, using different number of sensors $M$ for the IAEA 2D PWR case.}
    \label{fig: aFOM-ANL-OnlineTestError}
\end{figure}

To gain generality, an assessment of the reconstruction error in terms of $E_M$ and $\varepsilon_M$ must be performed. The results are shown in Figure \ref{fig: aFOM-ANL-OnlineTestError}: both TR-GEIM and PBDW have a good convergence rate. The thermal flux is the worst reconstructed field even though errors below 1\% are reached with less than 10 sensors.

\begin{figure}[htbp]
    \centering
    \includegraphics[width=0.9\linewidth]{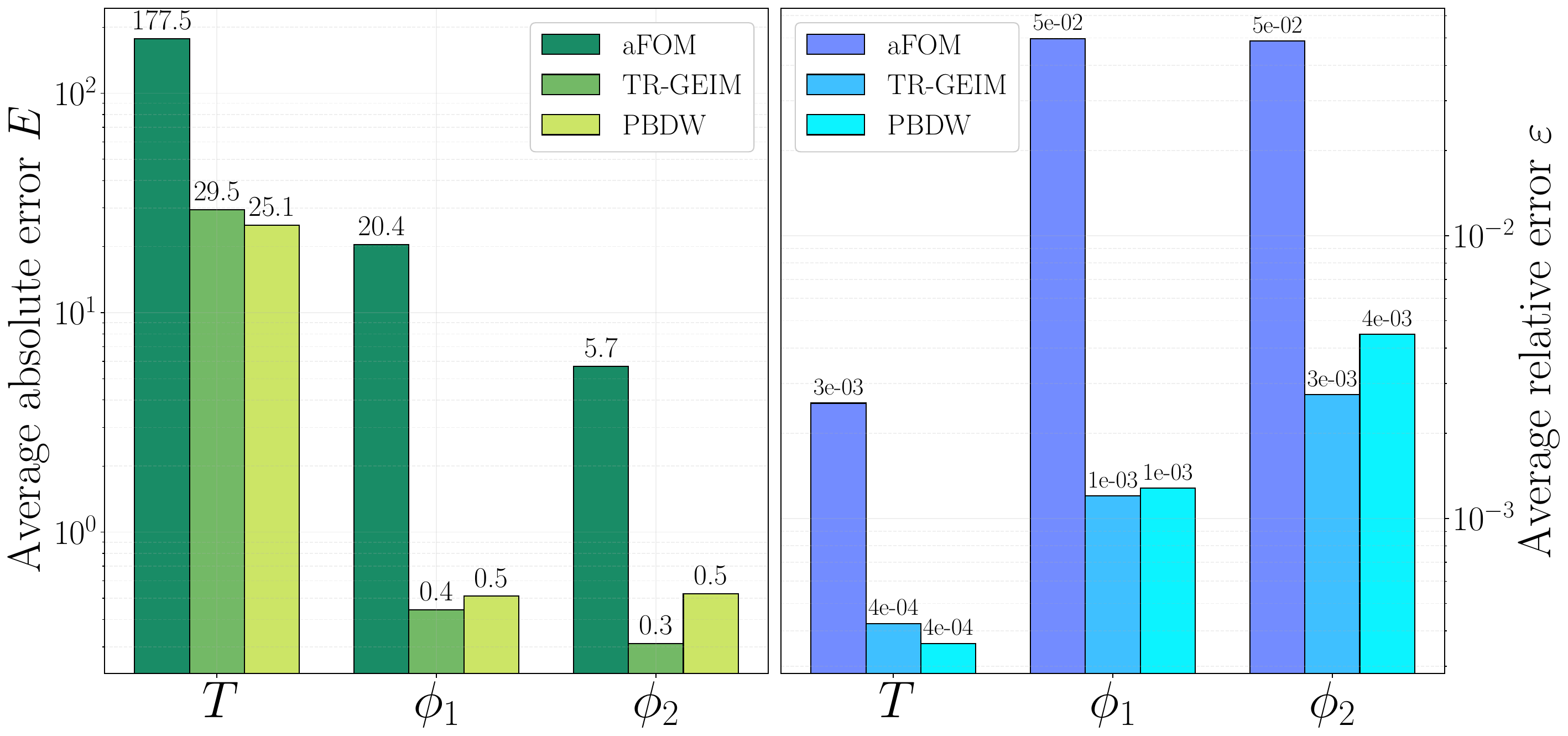}
    \caption{Bar plot (log-scale) of the absolute and relative error measured in $\norma{\cdot}_{L^2(\Omega)}$ using $M_{max}=15$ for the IAEA 2D PWR case.}
    \label{fig: aFOM-ANL-BarPlot}
\end{figure}

In the end, the absolute and relative errors of TR-GEIM and PBDW using 15 sensors are compared with the error of the aFOM compared to the \textit{ground-truth} (given by the FOM) in Figure \ref{fig: aFOM-ANL-BarPlot}. For every field, the improvement is important showing that the knowledge of the data is correctly combined within the reduced representation providing an update of the knowledge on the physical system. 

\subsubsection{TWIGL2D Benchmark}

The problem TWIGL2D-A described in Section \ref{sec-TWIGL2D-case-study} has the feedback effects only present in the absorption cross-section of group 1, Equation \eqref{eqn: TWIGL-A-feedback}. The parameter $\boldsymbol{\mu}$ is now an element of $\mathbb{R}^2$, including time $t$ and the feedback coefficients $\gamma_{a,1}$; in particular, $\boldsymbol{\mu}=[t,\gamma_{a,1}]^T\in\mathcal{D} = [0, 2]\times [10^{-3},10^{-2}]\subset\mathbb{R}^2$. The parameter space of time is divided into two sets: a training one $\Xi^{\text{train}}_t\subset(0,1]$ and a predict one $\Xi^{\text{predict}}_t\subset(1,2]$, sampled every 0.02 seconds in each subset; on the other hand, the training parameters of the feedback coefficient $\gamma_{a,1}$ are sampled in the whole interval, in particular, 20 values are uniformly selected, i.e. $\Xi^{\text{train}}_{\gamma_{a,1}}=[10^{-3}: 4.7\cdot 10^{-4}:10^{-2}]$, and 6 test parameters are sampled outside the training from, namely  $\Xi^{\text{predict}}_{\gamma_{a,1}}=[1\cdot 10^{-3} : 2.6\cdot 10^{-3}:1\cdot 10^{-2}]$. 

The number of training snapshots is $N_s = \text{dim}\left(\Xi^{\text{train}}_t\right)\cdot \text{dim}\left(\Xi^{\text{train}}_{\gamma_{a,1}}\right) = 50 \cdot 20 = 1000$. The offline phase has been performed to generate the basis functions and sensors; compared to the previous case, it has now been chosen to take $N=10$ as the dimension of the reduced space $Z_nN$ for the PBDW background space. Similarly to the previous case, the hyper-parameter $\xi$ has been tuned using a cross-validation procedure on a subset of the predict/test space.

\begin{figure}[htbp]
    \centering
    \includegraphics[width=0.98\linewidth]{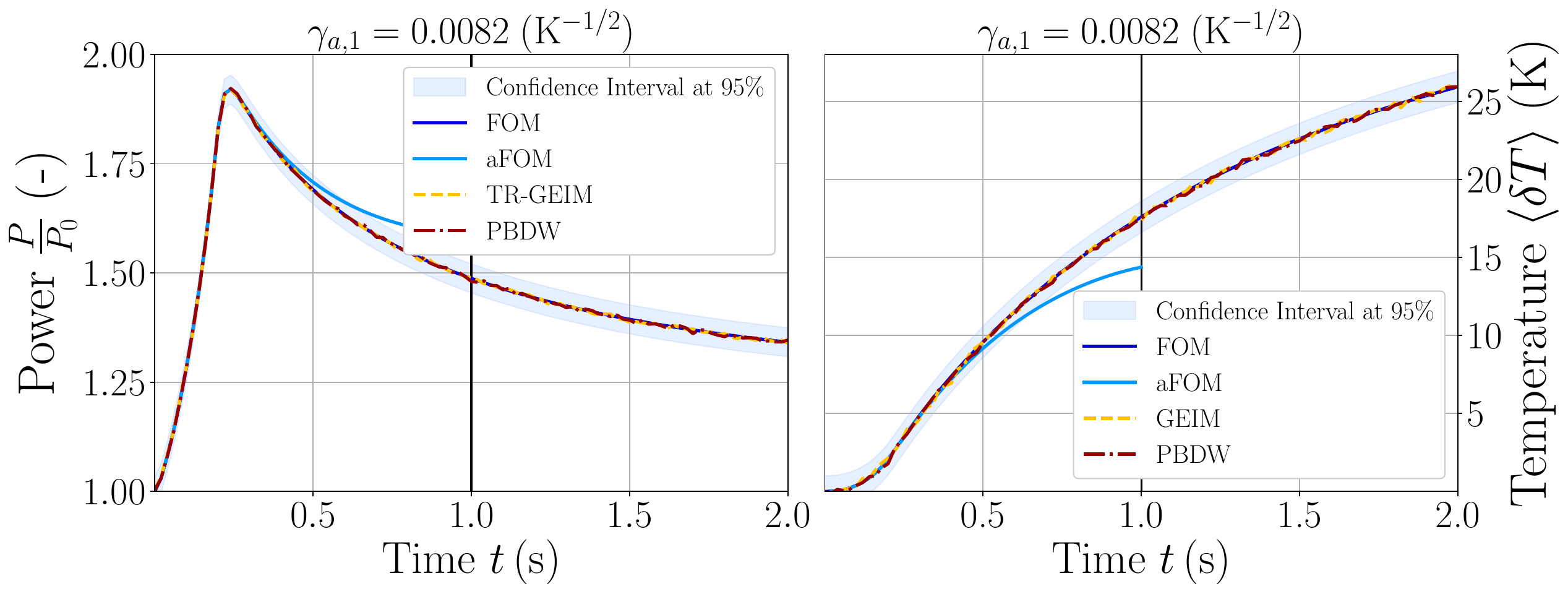}
    \caption{Line plots of the power (normalised with respect to the initial steady state) and the average temperature difference with respect to initial value for the TWIGL2D-A case, given $\gamma_{a,1} = 8.2\cdot 10^{-3}\,K^{-1}$.}
    \label{fig: aFOM-TWIGL-A-LinePlot}
\end{figure}

At first, the output of the data-driven ROM techniques, with $M_{max}=25$ measurements, is compared compared to the \textit{ground-truth} (FOM) and the training model (aFOM) in terms of total power and average temperature difference in Figure \ref{fig: aFOM-TWIGL-A-LinePlot} for $\gamma_{a,1} = 8.2\cdot 10^{-2}\,K^{-1}$. The inaccurate time evolution of the aFOM is corrected within TR-GEIM and PBDW and an accurate prediction behaviour is shown by both methods.

\begin{figure}[htbp]
    \centering
    \includegraphics[height=0.935\textheight]{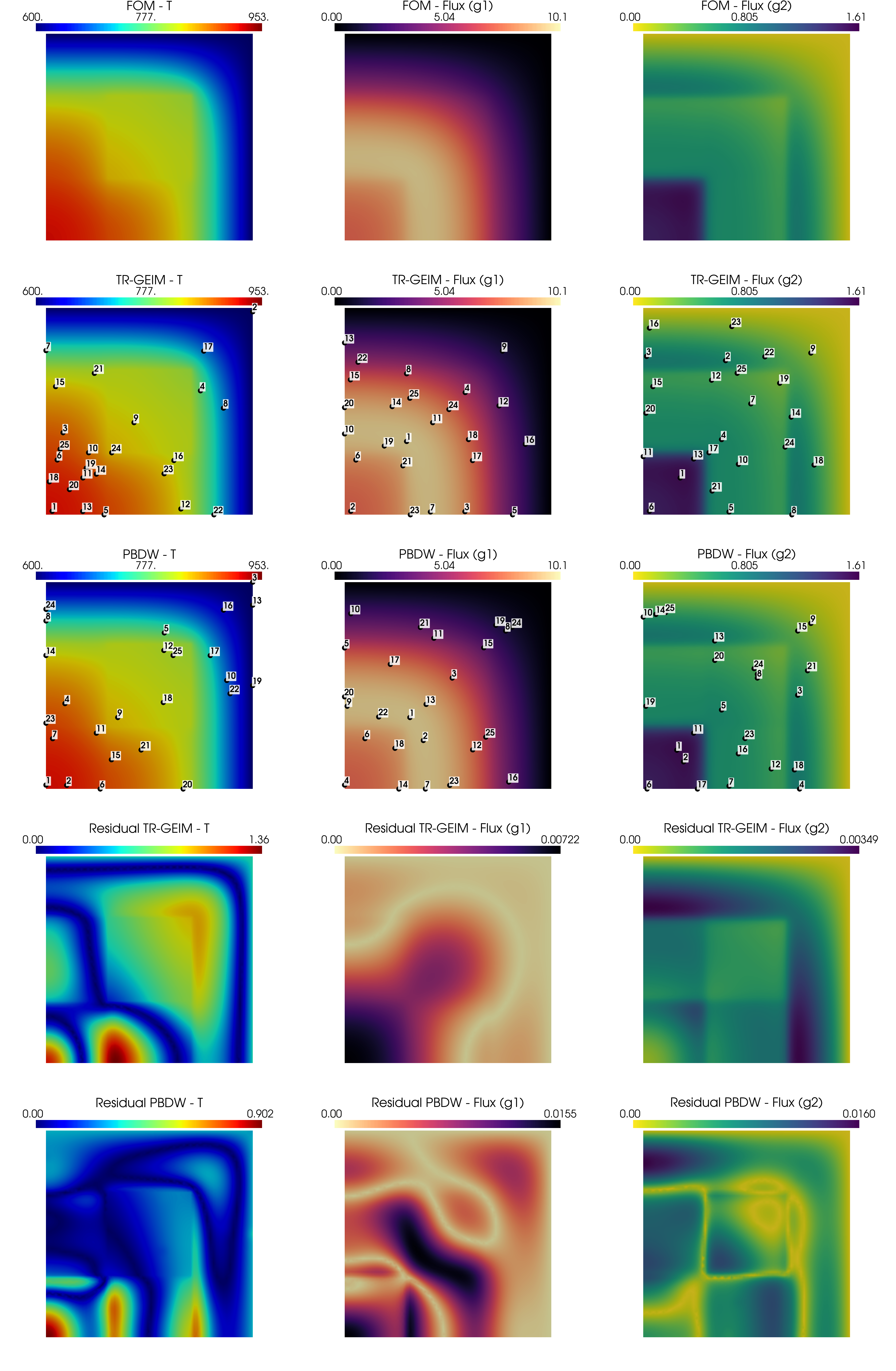}
    \caption{Contour plots of the temperature, fast flux and thermal flux fields (from left to right) for FOM, TR-GEIM and PBDW at final time $t= 2$ seconds, and the residual field for the TWIGL2D-A case for $\gamma_{a,1} = 8.2\cdot 10^{-3}\,K^{-1}$.}
    \label{fig: aFOM-TWIGL-A-Contour}
\end{figure}

The improvement can be also observed locally in Figure \ref{fig: aFOM-TWIGL-A-Contour} in which the reconstructions (using 25 measures) are plotted, along with the associated residual fields. The reconstruction of TR-GEIM and PBDW is very similar to the true field and the magnitude of the residual field is low compared to the value of the field itself; therefore, the information coming from the measurements is correctly combined with the background model, ensuring an improvement of the physical information.

\begin{figure}[htbp]
    \centering
    \includegraphics[width=1.\linewidth]{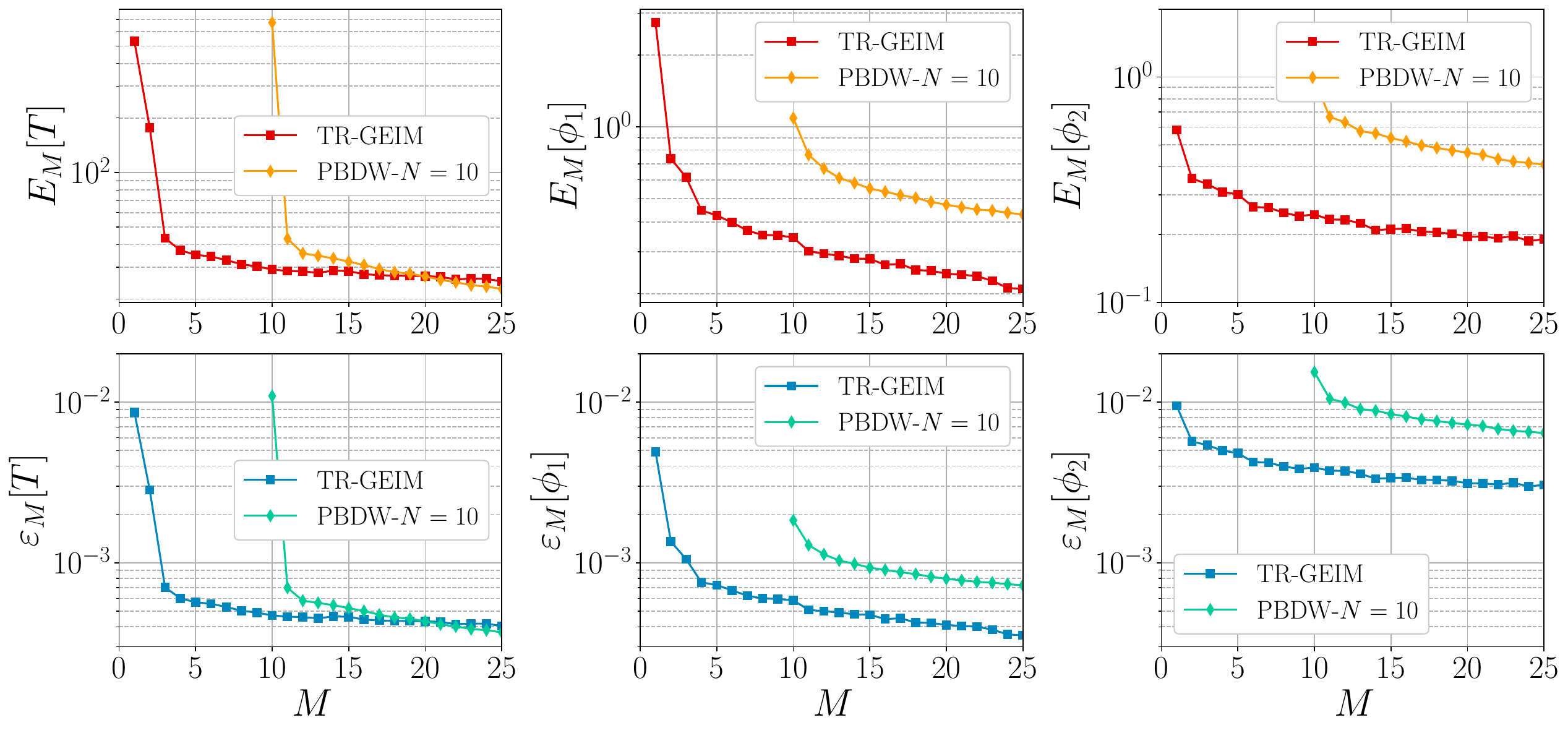}
    \caption{Average absolute and relative error measured in $\norma{\cdot}_{L^2(\Omega)}$, as in Eq. \eqref{eqn: error-defs}, using different number of sensors $M$ for the TWIGL2D-A case.}
    \label{fig: aFOM-TWIGL-A-OnlineTestError}
\end{figure}

Figure \ref{fig: aFOM-TWIGL-A-OnlineTestError} shows how the average absolute and relative error measured in $\norma{\cdot}_{L^2(\Omega)}$ behaves as the number of sensors increases: for all the fields the errors decrease in a good manner reaching sufficiently low values for the relative error. As before, the thermal flux seems to be more difficult to reconstruct compared to the others.

\begin{figure}[htbp]
    \centering
    \includegraphics[width=0.9\linewidth]{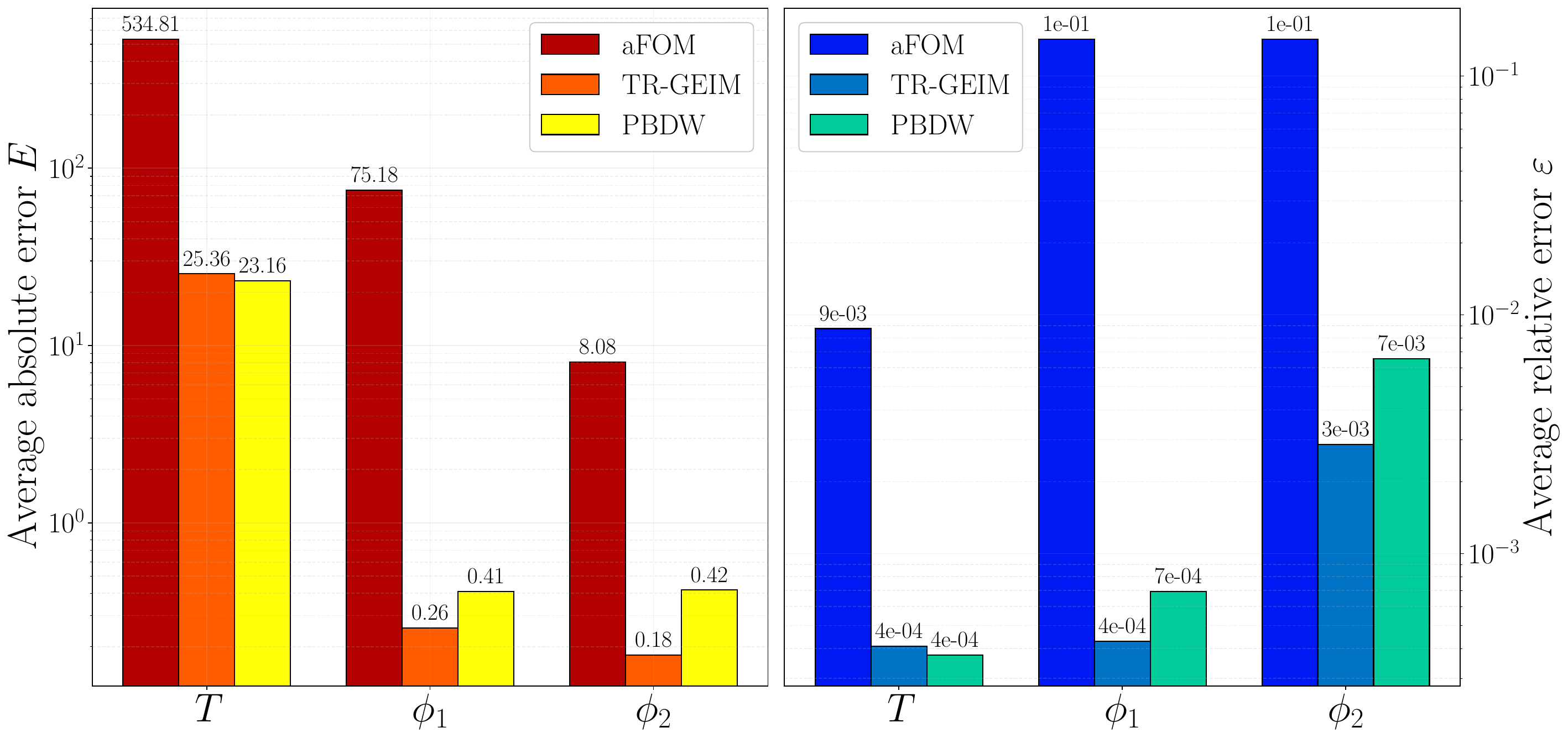}
    \caption{Bar plot (log-scale) of the average absolute and relative error measured in $\norma{\cdot}_{L^2(\Omega)}$ using $M_{max}=25$ for the TWIGL2D-A case.}
    \label{fig: aFOM-TWIGL-A-BarPlot}
\end{figure}

In the end, the error of the training model aFOM is compared with the one given by the TR-GEIM and PBDW using the maximum number of sensors available, i.e. $M_{max} = 25$. There is a significant improvement, especially for the neutron fluxes in which the error is greatly reduced by more than 2 orders of magnitude. 

\subsection{Model Correction from Linearised Coupling FOM}

\begin{figure}[htbp]
    \centering
    \includegraphics[width=0.9\linewidth]{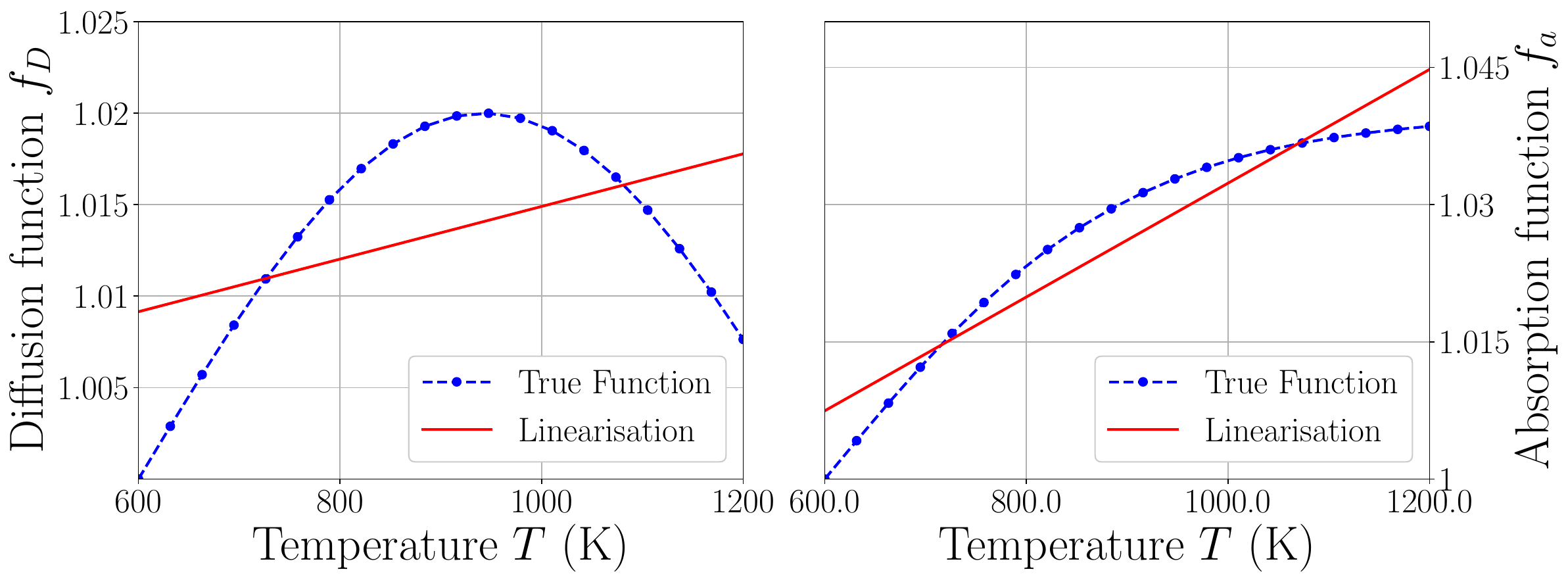}
    \caption{Coupling functions $f_D$ and $f_a$ and their approximation for LcFOM.}
    \label{fig:lin-coupling-LcFOM}
\end{figure}

Lastly, in this section, the TWIGL2D problem, described in Section \ref{sec-TWIGL2D-case-study}, is analysed from a different point of view: previously, the main focus was on the possibility of correcting the prediction of weakly coupled models, in which the single-physic codes were not able to directly communicate; whereas, in this last part the models are coupled together but the coupling functions are not known and hence they are approximated by linear relations, based on some training data-points. This case study will be referred to as TWIGL2D-B. In particular, the \textit{ground-truth} snapshots are generated using the coupling functions in Equation \eqref{eqn: TWIGL-B-feedback}, whereas, the snapshots of LcFOM are generated using a linear approximation of $f_D$ and $f_a$, as in Eq. \eqref{eqn: lin-coupling-fd-fa}, depicted in Figure \ref{fig:lin-coupling-LcFOM}: the one for the diffusion coefficients misses the real behaviour of the temperature dependence, the possibility of updating this information using data is investigated during the online phase.

The parameter $\boldsymbol{\mu}$ is an element of $\mathbb{R}^2$, including time $t$ and the diffusion coefficients $D_{1}$ in region 1; in particular, $\boldsymbol{\mu}=[t,D_1]\in\mathcal{D} = [0, 2]\times [0.5, 2]\subset\mathbb{R}^2$.  As before, the parameter space of time is divided into two sets: a training one $\Xi^{\text{train}}_t=(0,1]$ and a predict one $\Xi^{\text{predict}}_t=(1,2]$, sampled every 0.02 seconds. On the other hand, the training parameters of the diffusion coefficients $D_1$ are sampled in the whole interval, in particular 10 values are uniformly selected, i.e. $\Xi^{\text{train}}_{D_1}=[0.5: 0.166: 2]$; whereas, 5 test parameters are sampled outside the training from, namely  $\Xi^{\text{predict}}_{D_1}=[0.6 : 0.325:1.9]$. The number of training snapshots is $N_s = \text{dim}\left(\Xi^{\text{train}}_t\right)\cdot \text{dim}\left(\Xi^{\text{train}}_{\gamma_{a,1}}\right) = 51 \cdot 10 = 510$. After having performed the offline phase, the dimension of the reduced space $N$ of the PBDW background space is taken as 10. In the end, the hyper-parameter $\xi$ has been tuned using cross-validation as for the previous case studies.

\begin{figure}[htbp]
    \centering
    \includegraphics[width=0.98\linewidth]{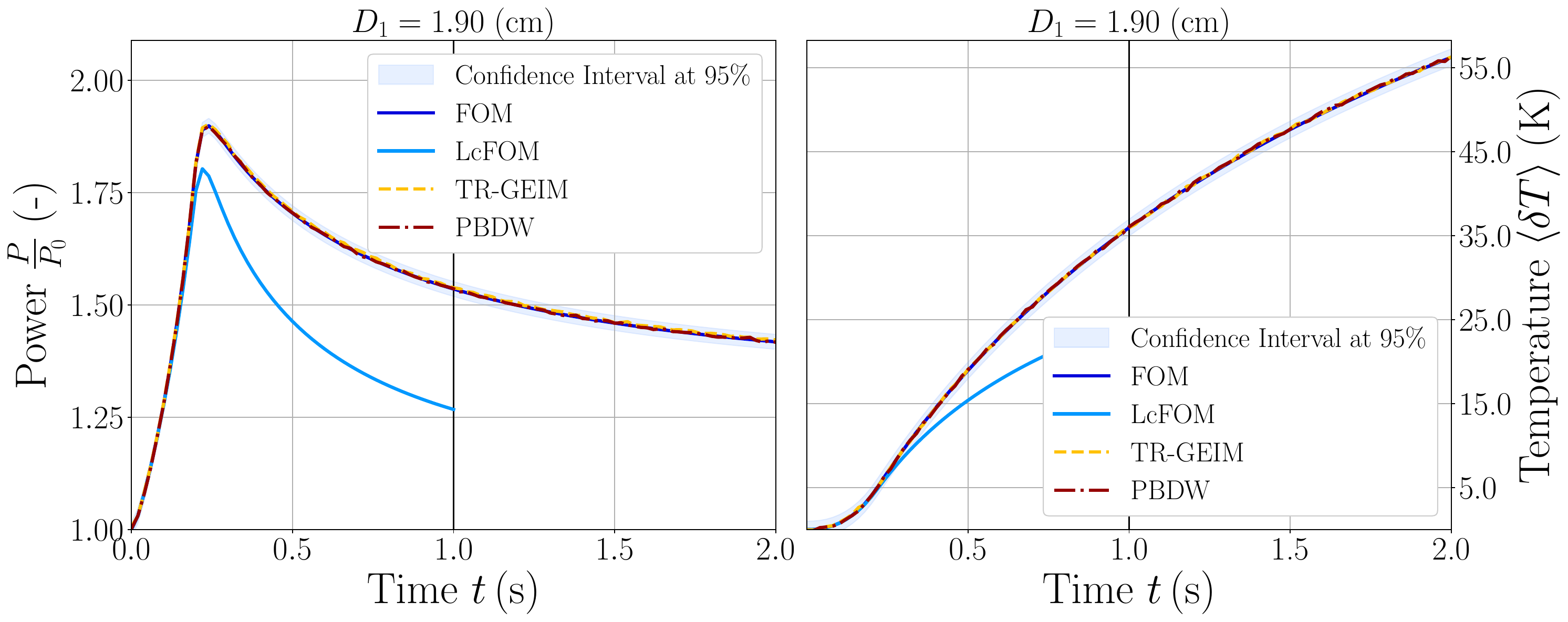}
    \caption{Line plots of the power (normalised with respect to the initial steady state) and the average temperature difference with respect to initial value for the TWIGL2D-B case with $D_1 = 1.9\,$cm.}
    \label{fig: LcFOM-TWIGL-B-LinePlot}
\end{figure}

Firstly, the power and the average temperature difference is plotted over time for the test value $D_1=1.9\,$cm in Figure \ref{fig: LcFOM-TWIGL-B-LinePlot} considering $M_{max}=25$ measurements. Compared to the previous section the LcFOM is much different in the evolution of the system, nevertheless, the TR-GEIM and the PBDW can retrieve the correct trend given by the FOM. Therefore, a wrong coupling function in the offline phase can be later corrected by the data during the online phase.

\begin{figure}[htbp]
    \centering
    \includegraphics[height=0.93\textheight]{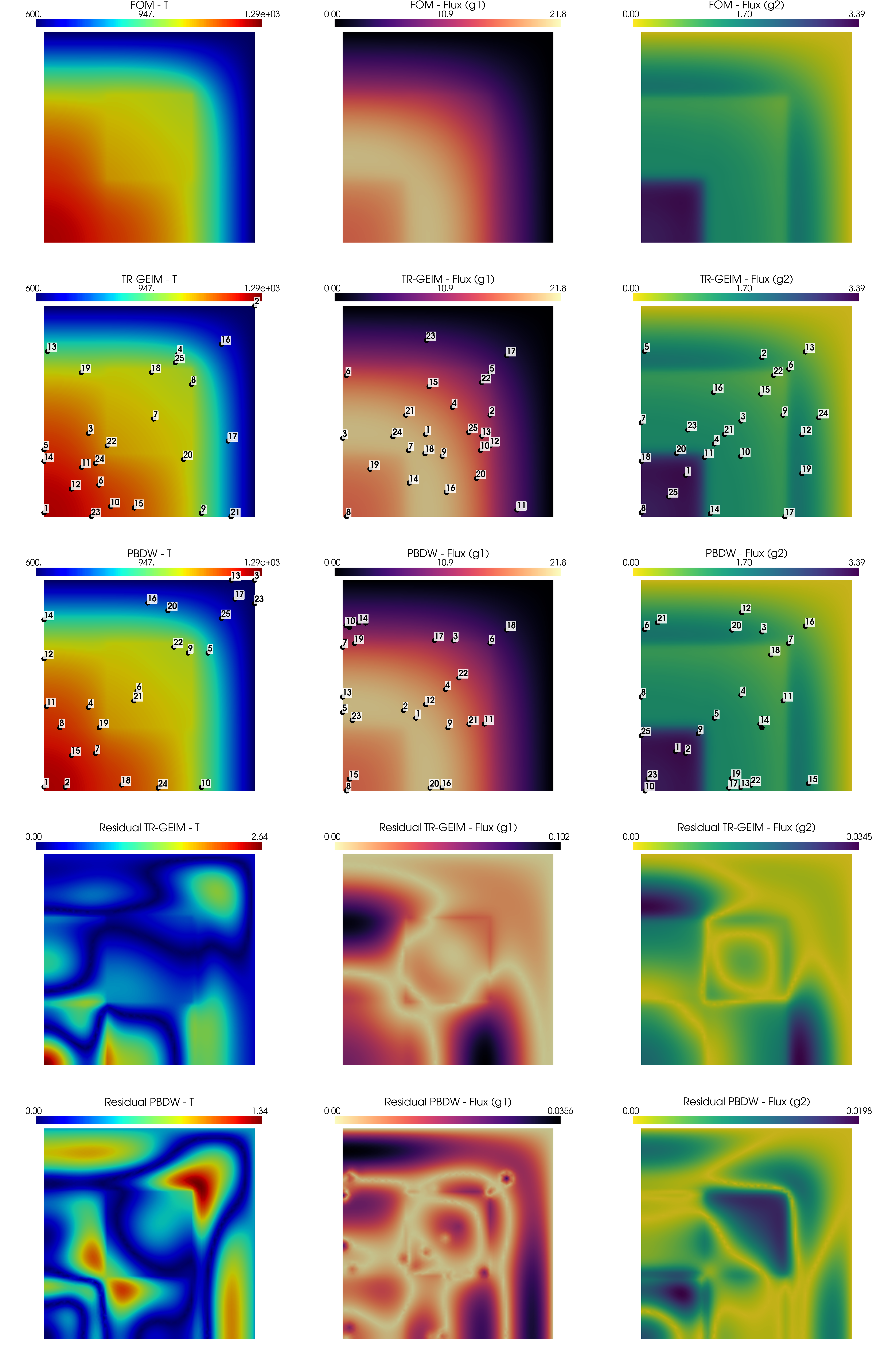}
    \caption{Contour plots of the temperature, fast flux and thermal flux fields (from left to right) for FOM, TR-GEIM and PBDW at final time $t$, and the residual field for the TWIGL2D-B case for $D_1 = 1.9\,$cm..}
    \label{fig: LcFOM-TWIGL-B-Contour}
\end{figure}

The update of the knowledge given by the data can be also observed by looking at the contour plots in Figure \ref{fig: LcFOM-TWIGL-B-LinePlot} at final time $t=2$ seconds for the same value of diffusion coefficient. The spatial dependence is very similar to the \textit{truth} and the related residual fields present low values, meaning that the reconstructed fields are very close to the FOM.

\begin{figure}[htbp]
    \centering
    \includegraphics[width=1.\linewidth]{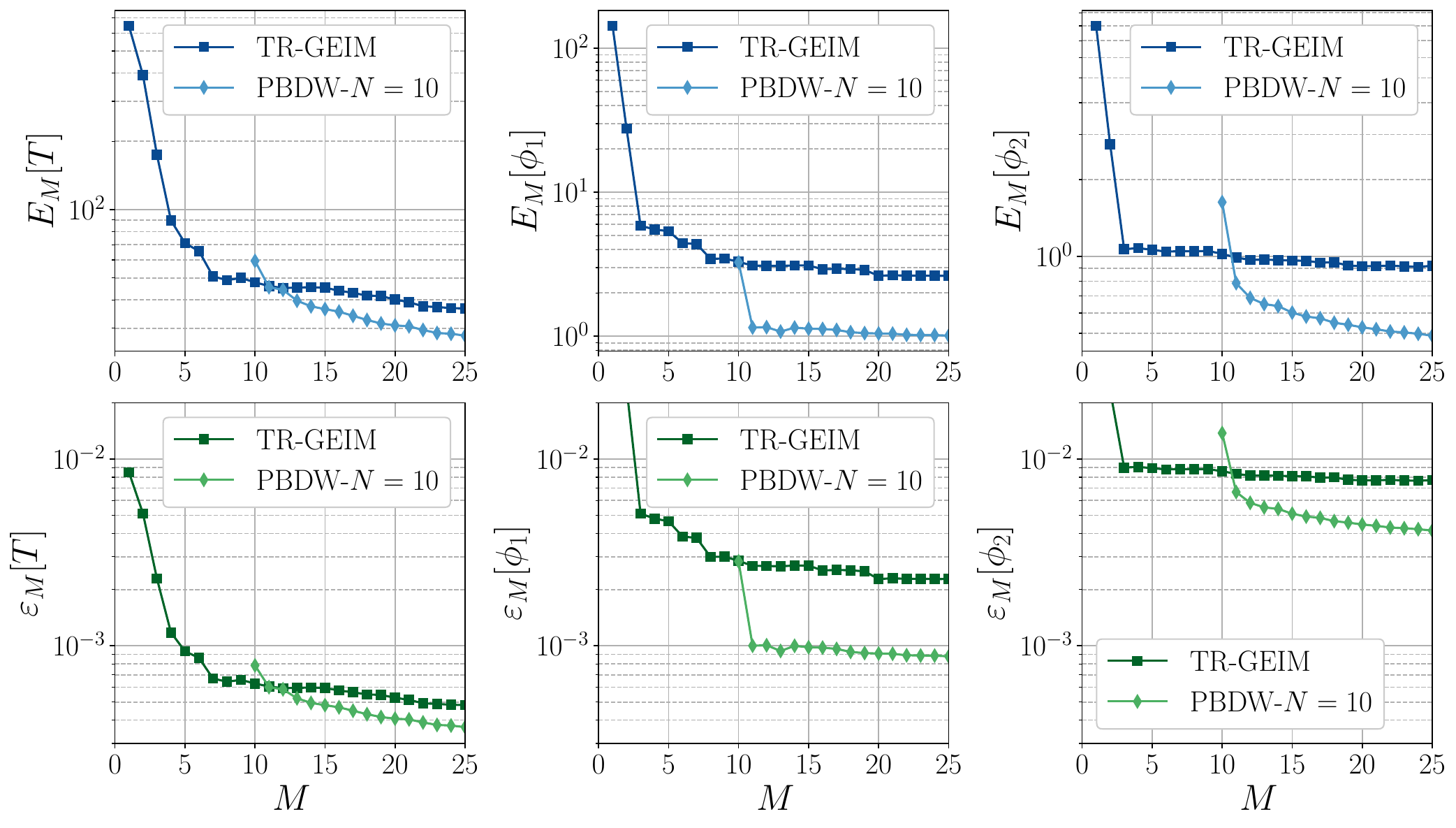}
    \caption{Average absolute and relative error measured in $\norma{\cdot}_{L^2(\Omega)}$, as in Eq. \eqref{eqn: error-defs}, using different number of sensors $M$ for the TWIGL2D-B case.}
    \label{fig: LcFOM-TWIGL-B-OnlineTestError}
\end{figure}

The numerical convergence of the methods can be observed in terms of the average absolute and relative error measured in $\norma{\cdot}_{L^2(\Omega)}$, reported in Figure \ref{fig: LcFOM-TWIGL-B-OnlineTestError}. The PBDW performs a bit better than the TR-GEIM, this behaviour is not shared by the analysis of the previous section: this can be related to a better approximation of the \textit{ground-truth} by the LcFOM, in particular, the reduced space built using POD is optimal in $L^2$-sense \cite{quarteroni2015reduced} compared to the one generated by GEIM. Hence the reduced space of the PBDW is closer to the truth, meaning that as more sensors are added the non-modelled physics is more efficiently introduced.

\begin{figure}[htbp]
    \centering
    \includegraphics[width=0.9\linewidth]{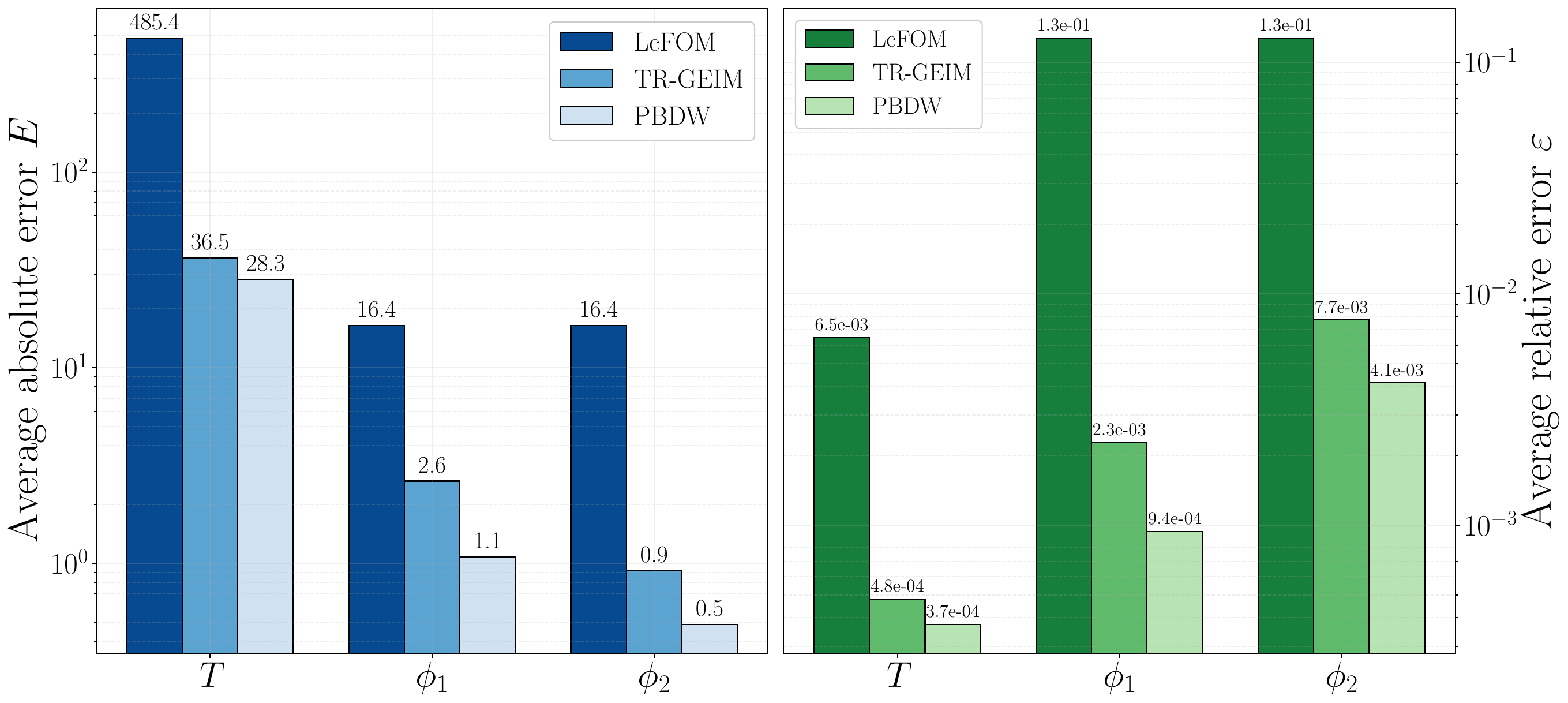}
    \caption{Bar plot (log-scale) of the absolute and relative error measured in $\norma{\cdot}_{L^2(\Omega)}$ using $M_{max}=25$ for the TWIGL2D-B case.}
    \label{fig: LcFOM-TWIGL-B-BarPlot}
\end{figure}

In the end, the error of TR-GEIM and PBDW is compared with the one performed by LcFOM considering 25 sensors and measurements per each test snapshot. The improvement is clearly visible showing that even a bad knowledge of the coupling function can be later corrected during the online phase.

\section{Conclusions}\label{sec: conclusions}

This work analyses the correction/update of model bias in the context of Multi-Physics modelling for nuclear applications combining the knowledge of noisy data (generated synthetically from the \textit{ground-truth}, i.e. a \textit{high-fidelity} fully coupled MP model) and the mathematical models subjected to simplifying assumptions, using the Tikhonov Regularised Generalised Empirical Interpolation Method and the Parameterised-Background Data-Weak formulation. 

These techniques have been applied to two case studies based on the IAEA 2D PWR and TWIGL2D benchmarks, analysing the model bias correction when single-physics codes are not fully coupled and when the coupling law is not known and it is approximated. In particular, different transients have been considered by selecting different reactivity insertions and, in both cases, the results obtained are really promising showing the strong capabilities of TR-GEIM and PBDW in updating the knowledge of the models and providing a corrected state estimation. 

In the future, further investigation of the outcomes of this work will be performed. In particular, these aspects will be scaled up to more realistic test cases including real data and considering other ways to introduce the coupling information, i.e. use of data-driven surrogate modelling.

\section*{Data Availability}
The \textbf{pyforce} package used in this work is part of the ROSE framework (Reduced Order multi-phySics data-drivEn) developed by the authors at the \href{https://github.com/ERMETE-Lab}{ERMETE Lab}. The package can be installed using \textit{pip} and will be soon available on Github at \href{https://github.com/ERMETE-Lab/ROSE-pyforce}{https://github.com/ERMETE-Lab/ROSE-pyforce} under the MIT license.

\section*{List of Symbols}
\textbf{Acronyms}

\begin{longtable}{rl}
\textbf{aFOM}    & approximated Full Order Model \\[5pt]
\textbf{DA}    & Data Assimilation \\[5pt]
\textbf{DDROM}    & Data Driven Reduced Order Modelling \\[5pt]
\textbf{FOM}    & Full Order Model \\[5pt]
\textbf{GEIM}    & Generalised Empirical Interpolation Method \\[5pt]
\textbf{IAEA}    & Internation Atomic Energy Agency \\[5pt]
\textbf{LcFOM}    & Linearised Coupling Full Order Model \\[5pt]
\textbf{MP}    & Multi-Physics \\[5pt]
\textbf{N-TH}    & Neutronics and Thermal-Hydraulics coupling\\[5pt]
\textbf{PBDW}    & Parameterised-Background Data-Weak formulation \\[5pt]
\textbf{PDE}    & Partial Differential Equation \\[5pt]
\textbf{POD}    & Proper Orthogonal Decomposition \\[5pt]
\textbf{POD-I}    & Proper Orthogonal Decomposition with Intepolation \\[5pt]
\textbf{PWR}    & Pressurised Water Reactor \\[5pt]
\textbf{ROM}    & Reduced Order Modelling \\[5pt]
\textbf{SVD}    & Singular Value Decomposition \\[5pt]
\textbf{TH}    & Thermal-Hydraulics \\[5pt]
\textbf{TR-GEIM}    & Tikhonov Regularised Generalised Empirical Interpolation Method \\[5pt]
\end{longtable}
\addtocounter{table}{-1}

\textbf{Greek Symbols}
\begin{longtable}{rl}
 {$\Omega$ } & {Spatial domain} \\
 {$\Sigma$ } & {Macroscopic cross section} \\
 {$\Upsilon$ } & {Library of sensors} \\
 {$\Xi$ } & {Subset of the parameter space} \\
 {$\alpha$ } & {PBDW coefficient of the background knowledge} \\
 {$\beta$ } & {GEIM reduced coefficient, \textit{inf-sup} constant, Delayed neutron fraction} \\
 {$\chi$ } & {Neutron fission spectrum} \\
 {$\delta$ } & {Tolerance for the GEIM greedy algorithm} \\
 {$\gamma$ } & {Parameters of the coupling function $f$} \\
 {$\epsilon$ } & {Random noise} \\
 {$\eta$ } & {Update in the PBDW} \\
 {$\lambda$ } & {TR-GEIM hyper-parameter, Precursors decay constant} \\
 {$\nu$ } & {Number of neutrons emitted by fission} \\
 {$\phi$ } & {Neutron flux} \\
 {$\rho$ } & {Density} \\
 {$\sigma$ } & {Noise standard deviation or generic standard deviation} \\
 {$\theta$ } & {PBDW coefficient of the update} \\
 {$\varepsilon$ } & {Relative error between the truth and a ROM reconstruction} \\
 {$\xi$ } & {PBDW hyper-parameter} \\
 {$\zeta$ } & {PBDW basis function} \\
 {$\boldsymbol{\mu}$ } & {Parameter vector} \\
\end{longtable}
\addtocounter{table}{-1}

\textbf{Latin symbols}
\begin{longtable}{rl}
 {$D$ } & {Diffusion coefficient} \\
 {$E$ } & {Absolute error between the truth and a ROM reconstruction} \\
 {$G$ } & {Number of neutronics energy groups} \\
 {$H$ } & {Heaviside step function} \\
 {$J$ } & {Number of neutronics precursors groups} \\
 {$L^p$ } & {Functional space of order $p$} \\
 {$M$ } & {Number of measures/sensors} \\
 {$N_s$ } & {Number of snapshots} \\
 {$P_0$ } & {Power normalisation constant} \\
 {$T$ } & {Temperature} \\
 {$U$ } & {Update space for PBDW} \\
 {$X$ } & {Hierarchical space built by GEIM} \\
 {$Z$ } & {Background space for PBDW} \\
 {$c$ } & {Neutron precursors} \\
 {$c_p$ } & {Specific heat capacity} \\
 {$f$ } & {Coupling function for N-TH} \\
 {$k$ } & {Thermal conductivity} \\
 {$g$ } & {PBDW basis sensor} \\
 {$q$ } & {GEIM magic function} \\
 {$q'''$ } & {Power density} \\
 {$r$ } & {Residual field} \\
 {$s$ } & {Point spread of the sensor} \\
 {$u$ } & {Generic Snapshot} \\
 {$v$ } & {Linear functional indicating a sensor, GEIM magic sensor, Neutron velocity} \\
 {$z$ } & {Background knowledge in the PBDW} \\
 {$\mathbb{A}$ } & {PBDW matrix of basis sensors} \\
 {$\mathbb{B}$ } & {GEIM matrix} \\
 {$\mathbb{C}$ } & {POD correlation matrix} \\
 {$\mathbb{E}$ } & {Expected value operator} \\
 {$\mathbb{K}$ } & {PBDW matrix of basis functions and sensors} \\
 {$\mathbb{I}$ } & {Identity matrix} \\
 {$\mathbb{P}$ } & {PBDW overall matrix} \\
 {$\mathbb{R}$ } & {Real number space} \\
 {$\mathbb{T}$ } & {Tikhonov GEIM matrix} \\
 {$\mathcal{D}$ } & {Parameter space} \\
 {$\mathcal{H}^p$ } & {Sobolev space of order $p$} \\
 {$\mathcal{K}$ } & {Kernel function} \\
 {$\mathcal{I}$ } & {GEIM interpolant} \\
 {$\mathcal{N}$ } & {Normal distribution} \\
 {$\mathcal{P}$ } & {Reconstruction operator} \\
 {$\mathcal{R}$ } & {Riesz representation operator} \\
 {$\mathcal{U}$ } & {Functional space of the snapshot} \\
 {$\vec{x}$ } & {Spatial coordinate} \\
 {$\vec{y}$ } & {Measurement vector} \\
\end{longtable}
\addtocounter{table}{-1}

\bibliography{bibliography.bib}
\clearpage

\appendix
\input{appendix}

\end{document}

%% file: schemes/mp_rom_start_art.tex
    \begin{tikzpicture}
        \node [draw,text centered,
            fill=bluePolimi!10,
            text width=2cm,
            minimum height=1.5cm
            ]  (FOM1) at (0,0) {FOM$_1$};
         
        \node [draw,text centered,
            fill=bluePolimi!10,
            text width=2cm,
            minimum height=1.5cm,
            right = 2cm of FOM1
            ]  (FOM2) {FOM$_2$};
         
        \node [draw,text centered,
            fill=greenPolimi!30,
            text width=2cm,
            minimum height=1.5cm
            ]  (ROM) at (2.2, -3.5) {MP ROM};
         
    \draw[->, line width = 1pt] (1.15,0.25) -- (3.15, 0.25)
    node[midway,right]{};
    \draw[->, line width = 1pt] (3.15,-0.25) -- (1.15, -0.25)
    node[midway,right]{};
    \node[text width=2cm] at (2.4,0.5) {Coupling};

    \draw[-{Triangle[width=18pt,length=8pt]}, line width=10pt](2.2,-1.5) -- (2.2, -2.75);
    
    \node[draw,bluePolimi,line width=2.5pt, inner sep=6.5mm,label=above:\color{bluePolimi}\textbf{Full Order MP model}\color{black},fit=(FOM1) (FOM2)] {};
    
\end{tikzpicture}

%% file: schemes/mp_rom_novel.tex
\begin{tikzpicture}
        \node [draw,text centered,
            fill=bluePolimi!10,
            text width=2cm,
            minimum height=1.5cm
            ]  (FOM1) at (0,0) {FOM$_1$};
         
        \node [draw,text centered,
            fill=bluePolimi!10,
            text width=2cm,
            minimum height=1.5cm,
            right = 2cm of FOM1
            ]  (FOM2) {FOM$_2$};
         
        \node [draw,text centered,
            fill=greenPolimi!30,
            text width=2cm,
            minimum height=1.5cm, 
            below = 2cm of FOM1
            ]  (ROM1) {ROM$_1$};
         
        \node [draw,text centered,
            fill=greenPolimi!30,
            text width=2cm,
            minimum height=1.5cm, 
            below = 2cm of FOM2
            ]  (ROM2) {ROM$_2$};
            
        \node [draw,text centered,
            fill=redPolimi!30,
            text width=2cm,
            minimum height=1.5cm, 
            below right = 1.cm and -0.1cm of ROM1
            ]  (Data) {Data};

    \draw[-{Triangle[width=15pt,length=8pt]}, line width=6pt](FOM1.south) -- (ROM1.north);
    \draw[-{Triangle[width=15pt,length=8pt]}, line width=6pt](FOM2.south) -- (ROM2.north);
    
    \draw[->, line width = 1pt] (1.15,-3.25) -- (3.15, -3.25)
    node[midway,right]{};
    \draw[->, line width = 1pt] (3.15,-3.75) -- (1.15, -3.75)
    node[midway,right]{};
    \node[text width=2cm] at (2.4,-4.05) {Coupling};
    
    \draw[->, dashed, line width = 1pt] (1.15,0.25) -- (3.15, 0.25)
    node[midway,right]{};
    \draw[->, dashed, line width = 1pt] (3.15,-0.25) -- (1.15, -0.25)
    node[midway,right]{};
    
    \draw[-{Triangle[width=15pt,length=8pt]}, line width=6pt, redPolimi](Data.north) -- (2.15, -4.2);
 
    \node[draw,greenPolimi,line width=2.5pt, inner sep=6.5mm,label=above:\color{greenPolimi}\textbf{MP-ROM}\color{black},fit=(ROM1) (ROM2)] {};
    
\end{tikzpicture}

%% file: schemes/pod-training.tex
\begin{tikzpicture}
        \node [draw,text centered,
            fill=greenPolimi!10,
            text width=4cm,
            minimum height=1.5cm
            ]  (N0) at (-3,0) {\textbf{Input Guess}: (N) solver with uniform $T$};
         
        \node [draw,text centered,
            fill=redPolimi!20,
            text width=4cm,
            minimum height=1.5cm,
            below = 2cm of N0
            ]  (POD-I-N) {POD-I on $\left\{\phi_g^{(i)}\right\}_{i=1}^{N_s}$};
         
        \node [draw,text centered,
            fill=greenPolimi!20,
            text width=4cm,
            minimum height=1.5cm,
            right = 2cm of POD-I-N
            ]  (TH) {(TH) solver with POD-I approximation of $q'''$};
         
        \node [draw,text centered,
            fill=redPolimi!20,
            text width=4cm,
            minimum height=1.5cm,
            below = 1.5cm of TH
            ]  (POD-I-TH) {POD-I on $\left\{T^{(i)}\right\}_{i=1}^{N_s}$};
         
        \node [draw,text centered,
            fill=greenPolimi!20,
            text width=4cm,
            minimum height=1.5cm,
            left = 2cm of POD-I-TH
            ]  (N) {(N) solver with POD-I approximation of $T$};
         
    \draw[->, draw = gray, dashed, line width = 1pt] (N0.south) -- (POD-I-N.north)
    node[near start,right]{Snapshots $\left\{\phi_g^{(i)}\right\}_{i=1}^{N_s}$};
    
    \draw[->, line width = 1pt] (POD-I-N.east) -- (TH.west)
    node[midway,below]{};
    \draw[->, dashed, line width = 1pt] (TH.south) -- (POD-I-TH.north)
    node[midway,left]{Snapshots $\left\{T^{(i)}\right\}_{i=1}^{N_s}$};
    \draw[->, line width = 1pt] (POD-I-TH.west) -- (N.east)
    node[midway,below]{};
    \draw[->, dashed, line width = 1pt] (N.north) -- (POD-I-N.south)
    node[midway,left] (arrow1) {Snapshots $\left\{\phi_g^{(i)}\right\}_{i=1}^{N_s}$};

    \node[draw,bluePolimi,line width=2.5pt, inner sep=6.5mm,label=above right:\color{bluePolimi}\textbf{aFOM solution loop}\color{black},fit=(POD-I-N) (POD-I-TH) (N) (TH) (arrow1)] {};
\end{tikzpicture}

%% file: appendix.tex
\section{Parameters of the case studies}

\subsection{IAEA 2D PWR benchmark}\label{app: pwr2d-params}
The values of the thermal heat capacity $c_p$ are not entirely physical, they have been chosen in order to have simple modelling of the system able to observe the feedback effects of the temperature on the neutronics in relatively slow time periods. The methods developed in this work are general and not bounded by the specific model adopted.
\begin{table}[htbp]
    \centering
    \begin{tabular}{c|cccc|c}
        \toprule
        Variable & Region 1 & Region 2 & Region 3 & Region 4 & Units\\ \midrule
        $D_1^{\text{ref}}\equiv D_1^{0}$ & 1.5 & 1.5 & 1.5 & 2 & (cm)\\
        $D_2^{\text{ref}}\equiv D_2^{0}$ & 0.4 & 0.4 &0.4 & 0.3 & (cm)\\
        $\Sigma_{a,1}^{0}$ & 0.01 & 0.01 & 0.01 & 0 & (cm) \\
        $\Sigma_{a,2}^{0}$ & 0.085 & 0.08 & 0.13 & 0.01 & (cm) \\
        $\Sigma_{s,1\rightarrow 2}$ & 0.02 & 0.02 & 0.02 & 0.04 & (cm)\\
        $\nu\Sigma_{f,1}$ & 0 & 0 & 0 & 0 & (n$\cdot$cm) \\
        $\nu\Sigma_{f,2}$ & 0.135 & 0.135 & 0.135 & 0 & (n$\cdot$cm) \\
        $B_{z,1}^2$ & $8\cdot10^{-5}$& $8\cdot10^{-5}$& $8\cdot10^{-5}$& $8\cdot10^{-5}$ & (cm$^{-2}$)\\
        $B_{z,2}^2$ & $8\cdot10^{-5}$& $8\cdot10^{-5}$& $8\cdot10^{-5}$& $8\cdot10^{-5}$ & (cm$^{-2}$)\\
        $\chi_1$ & 1   & 1   & 1 & 0 &(-) \\
        $\chi_2$ & 0   & 0   & 0 & 0 & (-) \\
        $v_1$ & $10^7$ & $10^7$ & $10^7$ & $10^7$ & (cm$/$s)\\
        $v_2$ & $10^5$ & $10^5$ & $10^5$ & $10^5$ & (cm$/$s)\\
        \midrule
        $k$ & 10.2 & 9 & 0.5 & 1.4 & $\left(\frac{\text{W}}{\text{cm K}}\right)$\\
        $\rho$ & 10.45 & 10.45 & 5. & 2.265 & $\left(\frac{\text{g}}{\text{cm}^3}\right)$\\
        $c_p$ & $235\cdot 10^{-6}$ & $235\cdot 10^{-6}$ & $100\cdot 10^{-6}$ & $1.4\cdot 10^{-6}$ & $\left(\frac{\text{J}}{\text{g K}}\right)$\\
        \bottomrule
    \end{tabular}
    \caption{Neutronic and thermal parameters for the IAEA 2D PWR problem.}
    \label{tab: anl-params-value}
\end{table}

\begin{table}[htbp]
    \centering
    \begin{tabular}{c|cc}
        \toprule
        Group $j$ & $\beta_j$ & $\lambda_j \;(\text{s}^{-1})$ \\ \midrule
        $1$ & 0.000247 & 0.0127 \\
        $2$ & 0.0013845 & 0.0317 \\
        $3$ & 0.001222 & 0.115 \\
        $4$ & 0.0026455 & 0.311 \\
        $5$ & 0.000832 & 1.4 \\
        $6$ & 0.000169 & 3.87 \\
        \bottomrule
    \end{tabular}
    \caption{Neutronic kinetics parameters for the IAEA 2D PWR problem.}
    \label{tab: anl-kin-params-value}
\end{table}

\subsection{TWIGL2D benchmark}\label{app: twigl-params}
The values of the thermal heat capacity $c_p$ are not entirely physical, they have been chosen in order to have simple modelling of the system able to observe the feedback effects of the temperature on the neutronics in relatively slow time periods. The methods developed in this work are general and not bounded by the specific model adopted.
\begin{table}[htbp]
    \centering
    \begin{tabular}{c|ccc|c}
        \toprule
        Variable & Region 1 & Region 2 & Region 3 & Units\\ \midrule
        $D_1^{\text{ref}}\equiv D_1^{0}$ & 1.4 & 1.4 & 1.3 & (cm)\\
        $D_2^{\text{ref}}\equiv D_2^{0}$ & 0.4 & 0.4 & 0.5 & (cm)\\
        $\Sigma_{a,1}^{\text{ref}}\equiv \Sigma_{a,1}^{0}$ & 0.01 & 0.01 & 0.008 & (cm) \\
        $\Sigma_{a,2}^{\text{ref}}$ & 0.15 & 0.15 & 0.05 & (cm) \\
        $\Sigma_{s,1\rightarrow 2}$ & 0.01 & 0.01 & 0.01 & (cm)\\
        $\nu\Sigma_{f,1}$ & 0.007 & 0.007 & 0.003 & (n$\cdot$cm) \\
        $\nu\Sigma_{f,2}$ & 0.2   & 0.2   & 0.06 & (n$\cdot$cm) \\
        $\chi_1$ & 1   & 1   & 1 & (-) \\
        $\chi_2$ & 0   & 0   & 0 & (-) \\
        $\beta$  & 0.0075 & 0.0075 & 0.0075 & (-)\\
        $\lambda$ & 0.08 & 0.08 & 0.08 & (s$^{-1}$)\\
        $v_1$ & $10^7$ & $10^7$ & $10^7$ & (cm$/$s)\\
        $v_2$ & $10^5$ & $10^5$ & $10^5$ & (cm$/$s)\\
        \midrule
        $k$ & 8 & 1 & 5 & $\left(\frac{\text{W}}{\text{cm K}}\right)$\\
        $\rho$ & 10.45 & 10.45 & 10.45 & $\left(\frac{\text{g}}{\text{cm}^3}\right)$\\
        $c_p$ & $235\cdot10^{-6}$ & $235\cdot10^{-6}$ & $235\cdot10^{-6}$ & $\left(\frac{\text{J}}{\text{g K}}\right)$\\
        \bottomrule
    \end{tabular}
    \caption{Neutronic and thermal parameters for the TWIGL2D problem.}
    \label{tab: twigl2d-params-value}
\end{table}

\section{Weak formulations and Finite Element for Multi-Group Neutron Diffusion}\label{app: FE-MG-diffusion}

The multi-group diffusion equations \eqref{eqn: mg-neutron-diffusion} , presented previously, is reported in the \textit{strong form}, meaning that the functions $\{\phi_g\}_{g=1}^G, \{c_j\}_{j=1}^J$ should be sufficiently regular, namely their second derivatives must exist and be continuous. Since it is not easy to prove that a solution of such regularity exists and is unique \cite{salsa}, the solution is searched in a broader functional space which imposes less regularity. Therefore, the \textit{strong} problem is rewritten in \textit{weak form}, so that the solution can be less regular\footnote{This procedure is consistent since, if there exists a solution to the strong problem, it satisfies the weak problem as well.}.

Before entering into the details of the \textit{weak} formulations, some basic tools are presented: at first, let $L^2(\Omega)$ be a Hilbert space s.t.
\begin{equation}
u\in L^2(\Omega) \Longleftrightarrow \int_\Omega \module{u}^2\,d\Omega <\infty;
\end{equation}
endowed with a scalar product
\begin{equation}
\scalarprod{u}{v}_{L^2} = \int_\Omega v\cdot u\,d\Omega \qquad\qquad \forall u,v\in L^2(\Omega).
\end{equation}
The \textit{weak} solutions are usually searched in the Sobolev functional spaces $\mathcal{H}^p(\Omega)$ (\cite{salsa, Gazzola_A3}), s.t.
\begin{equation}
\mathcal{H}^p(\Omega) = \left\{u\in L^2(\Omega)\,:\,\int_\Omega \module{D^{\bm{\alpha}}u}^2\,d\Omega <\infty \right\},
\end{equation}
given $D^{\bm{\alpha}}$ the multi-index derivative\footnote{Given $\bm{\alpha}\in\mathbb{N}^n$ of order $p=\sum_{i=1}^n \alpha_i$, the multi-index derivative of a function is defined as
\begin{equation*}
D^{\bm{\alpha}}u = \frac{\partial^pu}{\partial x_1^{\alpha_1} \dots \partial x_n^{\alpha_n}},
\end{equation*}
in which $n$ is usually 2 or 3. Moreover, it should be pointed out that the derivatives entering in the \textit{weak} formulation are meant to be \textit{weak} derivatives \cite{Gazzola_A3}.} of order $p$. Once this mathematical framework has been briefly presented, the \textit{weak} formulation can be obtained by testing the \textit{strong} form against a \textit{test function}, living in a suitable Sobolev space.

Let $\Gamma_D$ be the boundary onto which Dirichlet boundary conditions are imposed and let $\Gamma_N$ be the boundary onto which Neumann boundary conditions are imposed\footnote{Both homogeneous for the sake of simplicity.}, so that $\Gamma_D\cup \Gamma_N = \partial \Omega$; the neutron flux $\phi_g\in\mathcal{H}^1(\Omega)$ with $\left.\phi_g\right|_{\Gamma_D}=0$ satisfies the weak diffusion equation, given $\psi_g\in\mathcal{H}^1(\Omega)$ with $\left.\psi_g\right|_{\Gamma_D}=0$ as the test function
\begin{equation}
    \begin{split}
        \frac{1}{v_g}\int_\Omega\dpart{\phi_g}{t}\cdot \psi_g\,d\mathbf{x} &+ \int_\Omega D_g \nabla \phi_g\cdot \nabla \psi_g\,d\mathbf{x}\\
        &+\int_\Omega \left(\Sigma_{a,g}
        +\sum_{g'\neq g}\Sigma_{s,g\rightarrow g'}+D_gB_{z,g}^2\right)\phi_g\cdot \psi_g\,d\mathbf{x}\\
        &-\sum_{g'\neq g}\int_\Omega \Sigma_{s,g'\rightarrow g}\phi_{g'}\cdot \psi_g\,d\mathbf{x}\\
        &-\chi_g^p\frac{1-\beta}{k_{\text{eff}}}\sum_{g'}\int_\Omega \nu_{g'}\Sigma_{f,g'}\phi_{g'}\cdot \psi_g\,d\mathbf{x}\\
        &-\sum_{j=1}^J\chi_g^{d_j}\int_\Omega \lambda_jc_j\cdot \psi_g\,d\mathbf{x}=0
    \end{split}
    \label{eqn: weak-trans-diff}
\end{equation}
whereas the precursors $c_j\in L^2(V)$ satisfy the weak precursors equation given $\eta_j\in L^2(V)$ as the test function 
\begin{equation}
    \int_\Omega\dpart{c_j}{t}\cdot \eta_j\,d\mathbf{x} -
    \frac{\beta_j}{k_{\text{eff}}}\sum_{g}\int_\Omega\nu_g\Sigma_{f,g}\phi_{g}\cdot \eta_j\,d\mathbf{x} + \int_\Omega \lambda_jc_j\cdot \eta_j\,d\mathbf{x}=0
\end{equation}